\documentclass[ejs]{imsart}

\usepackage[pdftex]{graphicx} %Loading the package
\graphicspath{{}{figures/}}

\usepackage{amsthm,amsmath,amssymb}
\usepackage[normalem]{ulem}
\usepackage{bm,bbm}

\usepackage{caption,subcaption}
\usepackage{multirow,booktabs,array}\usepackage[color,notref,notcite,dvipsnames]{xcolor}
\definecolor{labelkey}{rgb}{0.6,0,1}

\usepackage{natbib}
\usepackage{appendix}

\newtheorem{theorem}{Theorem}[section]
\newtheorem{lemma}{Lemma}[section]
\newtheorem{proposition}{Proposition}[section]
\newtheorem{definition}{Definition}[section]

\newtheorem{example}{Example}[section]
\newtheorem{remark}{Remark}[section]
\newtheorem{corollary}{Corollary}[section]
\newtheorem{conjecture}{Conjecture}[section]
\newtheorem{assumption}{Assumption}[section]
\newtheorem{algorithm}{Algorithm}[section]

\numberwithin{equation}{section}
\numberwithin{table}{section}
\numberwithin{figure}{section}

\DeclareMathOperator{\E}{E}

\startlocaldefs

\endlocaldefs

\def\iotaa{\bm{\iota}}

\begin{document}

\begin{frontmatter}
% "Title of the paper"
\title{Clustered Archimax copulas}
\runtitle{Clustered Archimax copulas}
\begin{aug}
	\author{\fnms{Simon} \snm{Chatelain}}\thanksmark{t1}\thanksmark{t2}\thanksref{r1}
	\quad
	\author{\fnms{Samuel} \snm{Perreault}}\thanksmark{t3}
	\\
	\author{\fnms{Anne-Laure} \snm{Foug\`eres}}\thanksmark{t2}
	\quad
	\author{\fnms{Johanna} \snm{G. Ne{\v{s}}lehov\'a}}\thanksmark{t4}
	\runauthor{Chatelain et al.}
	
	\thankstext{r1}{Corresponding author: simon.chatelain.ev@gmail.com}	
		
	\address{
	\thanksmark{t1}Geneva School of Economics and Management, Universit\'e de Gen\`{e}ve\\
	\thanksref{t2}Universit\'{e} de Lyon, Universit\'{e} Lyon 1, CNRS UMR 5208, Institut Camille-Jordan\\
	\thanksref{t3}Department of Statistical Sciences, University of Toronto\\
	\thanksref{t4}Department of Mathematics and Statistics, McGill University}
\end{aug}

\begin{abstract}
When modeling multivariate phenomena, properly capturing the joint extremal behavior is often one of the many concerns.
Archimax copulas appear as successful candidates in case of asymptotic dependence.
In this paper, the class of Archimax copulas is extended via their stochastic representation to a clustered construction.
These clustered Archimax copulas are characterized by a partition of the random variables into groups linked by a radial copula; each cluster is Archimax and therefore defined by its own Archimedean generator and stable tail dependence function.
The proposed extension allows for both asymptotic dependence and independence between the clusters, a property which is sought, for example, in applications in environmental sciences and finance.
The model also inherits from the ability of Archimax copulas to capture dependence between variables at pre-extreme levels.
The asymptotic behavior of the model is established, leading to a rich class of stable tail dependence functions.
\end{abstract}

%\begin{keyword}[class=MSC]
%\kwd[Primary ]{62H12} 
%\kwd{62G05}
%\kwd{62G20}
%\kwd{62G32}
%\kwd[; secondary ]{60G70}
%\end{keyword}


\end{frontmatter}

\thispagestyle{empty}
% AOS,AOAS: If there are supplements please fill:

\section{Introduction}\label{sec:intro}
The Archimax copula model introduced by \cite{Caperaa/Fougeres/Genest:1997} has been advocated in \cite{Chatelain/Fougeres/Neslehova:2020} as a flexible way to model  a group of variables whose asymptotic dependence is driven by a stable tail dependence function (stdf); or more precisely a random vector whose dependence structure follows an asymptotic extreme-value regime perturbed by the same distortion. However, as it is the case for rainfall over large territories for example, asymptotic independence between certain variables is likely to be present {and this phenomenon cannot be handled by a single Archimax model without limiting the marginal dependence structure to be an (exchangeable) Archimedean copula. In financial applications, a stock portfolio might contain stocks from the same industry, causing them to be dependent in the extreme regime, while stocks from different industries might be asymptotically independent. Likewise, assuming the same distortion for all variables may not be realistic when the number of variables is large.
	
We propose a dependence model based on the survival copula of a nonnegative, real random vector of the form
\begin{equation*} %\label{eq:stochrep}
	\bm{X}=({R}_1\bm{S}_1,\ldots,R_K\bm{S}_K)\;,
\end{equation*}
where for all $k\in\{1,\ldots,K \}$, $R_k$ are nonnegative random variables and $\bm{S}_k$ are $d_k$-dimensional random vectors which will be characterized later. More specifically, each of the vectors $R_k\bm{S}_K = R_k\times(S_{k1},\ldots,S_{k d_k})$ will be defined so as to have an Archimax survival copula, while conditional independence of the vectors $\{\bm{S}_k\}_{k=1}^K$ given $\{R_k\}_{k=1}^K$ will introduce parsimony to the dependence model. This paper focuses on the survival copula of $\bm{X}$, which can subsequently be used in a copula model with arbitrary margins.

Stochastic representations involving a randomly scaled random vector are very common in the literature. For example, pseudo-polar representations exist for elliptical (see \cite{Fang/Kotz/Ng:1990}), generalized Pareto (see \cite{Ferreira/deHaan:2014}), Archimedean (see \cite{McNeil/Neslehova:2009}) and Liouville (see \cite{Belzile/Neslehova:2017}). Such constructions have garnered significant interest as models suited for extreme values; they can offer the flexibility needed to model across extremal dependence classes as shown in \cite{Huser/Opitz/Thibaud:2017}, \cite{Wadsworth/Tawn/Davison/Elton:2017} and \cite{Huser/Wadsworth:2019}. \cite{Engelke/Opitz/Wadsworth:2019} offer an extensive study of the extremal behavior of bivariate pseudo-polar vectors. As shown in \cite{Charpentier/Fougeres/Genest/Neslehova:2014}, Archimax copulas also allow a representation of the form $R\bm{S}$; the interpretation proposed by \cite{Chatelain/Fougeres/Neslehova:2020} is that of a radial variable $R$ distorting the vector $\bm{S}$ characterized by a stdf and therefore representing the extremal regime. The aim of this paper is to propose a dependence model in a way that its higher-dimensional margins are  Archimax copulas but with possibly different distortions or stdfs. There is also a connection to be made with \cite{Hofert/Huser/Prasad:2017} which have extended the class of Archimax copulas with completely monotone Archimedean generators to hierarchical constructions. Therein, hierarchies can be introduced either via the frailties or stdfs, but the extremal behavior is not elicited.

In this paper, we define the new class of clustered Archimax copulas. Building a model from this family is straightforward as it only requires specifying the Archimax clusters and the dependence between the radial variables $R_1,\ldots,R_K$. Under non-restrictive assumptions, the extremal behavior of the clustered Archimax copula is obtained and shown to be quite flexible. We find that different choices of popular Archimedean generators coupled with dependence models on the radial variables lead to various regimes falling under either asymptotic dependence or independence between clusters. Meanwhile, each Archimax cluster retains its flexibility in the extreme dependence regime. We also propose inference techniques to fit clustered Archimax copulas to data and illustrate with a data application on rainfall in France.

The paper is organized as follows. Section~\ref{sec:prelim} contains preliminary notions; Section~\ref{sec:modelproperties} defines the model;  Section~\ref{sec:modelproperties} presents results on its extremal behavior; Section~\ref{sec:simulations} provides illustrative examples; Section~\ref{sec:inference} covers inference techniques; Section~\ref{sec:data-application} contains the data application; finally, Section~\ref{sec:discussion} concludes the paper with a discussion. Proofs are reported in Appendix~\ref{sec:proofsofmodelproperties}, while Appendix~\ref{sec:conjecture} formulates a conjecture extending Theorem~\ref{thm:MaxDofA}. 

\section{Preliminaries}\label{sec:prelim}
Copulas contain all the information pertaining to the dependence between the components of continuous random vectors. The decomposition of \cite{Sklar:1959} links the marginals and the joint distribution of a random random vector via a copula, which is simply a distribution function on the unit hypercube with standard uniform margins. Consider $\bm{X}=(X_1,\ldots,X_d)\sim F$ with margins $F_1,\ldots,F_d$. Then there exists a copula $C$ such that for any $(x_1,\ldots,x_d)\in\mathbb{R}^d$,
$$
F(x_1,\ldots,x_d)=\Pr(X_1\le x_1,\ldots,X_d\le x_d)=C\{F_1(x_1),\ldots,F_d(x_d)\}\;.
$$
Moreover, $C$ is unique if the margins are all continuous. There is also a version of Sklar's decomposition for survival functions. That is, given the marginal survival functions $\bar{F}_1,\ldots,\bar{F}_d$ and the joint marginal function $\bar{F}$, there exists a copula $\bar{C}$, called a survival copula of $\bm{X}$, such that for all $(x_1,\ldots,x_d)\in\mathbb{R}^d$,
$$
\bar{F}(x_1,\ldots,x_d)=\Pr(X_1> x_1,\ldots,X_d> x_d)=\bar{C}\{\bar{F}_1(x_1),\ldots,\bar{F}_d(x_d)\}\;.
$$
It can be shown that if $\bm{U}\sim C$, then $1-\bm{U}\sim \bar{C}$ (see \cite{Nelsen:2006}). There is a large amount of literature on dependence modeling via copulas; we refer to the monographs of \cite{Nelsen:2006}, \cite{Durante/Sempi:2010} and \cite{Joe:2015}. 

One particular family of copulas, called Archimax, was introduced in two dimensions by \cite{Caperaa/Fougeres/Genest:2000}, and extended to arbitrary dimensions by \cite{Mesiar/Jagr:2013} and \cite{Charpentier/Fougeres/Genest/Neslehova:2014}. {Generalizing both Archimedean and extreme-value copula families, Archimax copulas take the form, for all $(u_1,\ldots,u_d)\in(0,1)$},
\begin{equation}\label{eq:axc}
C_{\psi,\ell}(u_1,\ldots,u_d)=\psi\circ \ell\{\phi(u_1),\ldots,\phi(u_d) \}\;,
\end{equation}
where $\psi$ is an Archimedean generator with inverse $\phi$ and $\ell$ is a stdf. These two functional parameters are defined below.
\begin{definition}\label{def:axc}
	A non-increasing and continuous function $\psi: [0,\infty) \to [0,1]$ which satisfies  $\psi(0)=1$, $\lim_{x\to\infty}\psi(x)=0$ and is strictly decreasing on $[0,x_\psi)$, where $x_\psi =\inf\{x:\psi(x)=0\}$, is  called an  Ar\-chi\-me\-de\-an generator. 
	
	A function $\ell : \mathbb{R}_+^d \to \mathbb{R}_+$ is called a $d$-variate stable tail dependence function (stdf) if there exists a finite measure $\mu$ on the $d$-dimensional unit simplex $\Delta_d=\{\bm{w}\in[0,1]^d:w_1+\cdots+w_d=1\}$ satisfying $\int_{\Delta_d}s_j d\mu(\bm{s}) = 1$ for all $j \in \{1,\ldots, d\}$, such that for all $\bm{x} \in \mathbb{R}^d_+$, 
	$$
	\ell(\bm{x}) = \int_{\Delta_d} \max(x_1s_1,\dots, x_d s_d) d\mu(\bm{s}).
	$$
	
A $d$-dimensional copula $C$ is called Ar\-chi\-max  if it permits the representation \eqref{eq:axc} for some $d$-variate stdf $\ell$ and an Ar\-chi\-me\-de\-an generator $\psi$ with inverse $\phi:(0,1]\to [0,\infty)$, where by convention $\psi(\infty)=0$ and $\phi(0)=x_\psi$. 
\end{definition}
In the special case when $\psi(x)=\psi_\Pi(x)=\exp(-x)$, one has that $C_{\psi,\ell}$ reduces to the extreme value copula with stdf $\ell$ defined for all $(u_1,\ldots,u_d)\in(0,1)^d$ by
$$
C_{\ell}(\bm{u})=\exp\{-\ell(-\ln u_1,\ldots,-\ln u_d)  \}\;.
$$
Stable tail dependence functions were introduced by \cite{Huang:1992} and are given a characterization by \cite{Ressel:2013} in terms of homogeneity, convexity and boundary properties. In this paper, the $d$-norm representation of stdfs is particularly useful. The following characterization, as discussed in \cite{Aulbach/Falk/Zott:2015}, can be traced back to the work of \cite{Pickands:1975}, \cite{deHaan/Resnick:1977}, and \cite{Vatan:1985} on the representation of standard max-stable processes. Any $d$-dimensional stdf $\ell$ can be written, for $(x_1,\ldots,x_d)\in \mathbb{R}_+^d$, as
\begin{equation}\label{eq:dnorm}
\ell(x_1,\ldots,x_d)= \E\bigl(\max_{1 \leq k \leq d} x_kW_k\bigr)
\end{equation}
for some positive random variables $W_1,\ldots,W_d$ with unit mean. When $\ell(\bm{x})=\ell_\Pi(\bm{x})=x_1+\ldots+x_d$, i.e., the stdf corresponding to asymptotic independence, then $C_{\psi,\ell}$ is simply an Archimedean copula given for all $(u_1,\ldots,u_d)\in(0,1)^d$ by $$C_{\psi}(\bm{u})=\psi\{\phi(u_1)+\ldots+\phi(u_d) \}\;.$$  Conditions on $\psi$ for $C_\psi$ to be a copula were explored in \cite{McNeil/Neslehova:2009}, while conditions for $\psi$ and $\ell$ for \eqref{eq:axc} to be a copula were explored in \cite{Charpentier/Fougeres/Genest/Neslehova:2014}.

Archimax copulas also admit a stochastic representation, which this paper builds upon. Theorem~3.3 of \cite{Charpentier/Fougeres/Genest/Neslehova:2014} states under conditions on $\psi$ and $\ell$, that $C_{\psi,\ell}$ is the survival copula of a random vector 
\begin{equation}\label{eq:rep}
(X_1,\ldots,X_d)= (R S_1,\ldots,R S_d) = R \bm{S},
\end{equation}
where $R$ is a positive random variable independent of $\bm{S} = (S_1,\dots,S_d)$. The distribution of $R$ is linked to the Archimedean generator via the Williamson-$d$ transform, i.e. $\psi=\mathfrak{W}_d (F_R)$ and $F_R=\mathfrak{W}_d^{-1}(\psi)$. We refer to \cite{McNeil/Neslehova:2009} and \cite{Larsson/Neslehova:2011} for more details. The survival function of $\bm{S}$ is given, for any $\bm{s}\in\mathbb{R}_+^d$, by
\begin{equation}\label{eq:Ssurvivalfunction}
\bar{G}_\ell(\bm{s})=\Pr(S_1>s_1,\ldots,S_d>s_d)=\left[\max\{0,1-\ell(\bm{s})\}\right]^{d-1}\;.
\end{equation}
Note in particular that the margins of $\bm{S}$ are Beta. Specifically, $S_i \sim \mathrm{Beta}(1,d-1)$ for all $i \in  \{1,\ldots, d\}$. We interpret \eqref{eq:rep} as a dependence structure defined by a distortion (or radial) random variable $R$ applied to the extremal component $\bm{S}$.

Archimax copulas have a given extreme-value attractor, which motivates their use to model pre-extreme dependence. Recall that a function $f: \mathbb{R}_+ \to \mathbb{R}_+$ is regularly varying with index $\alpha \in \mathbb{R}$ if and only if for all $x > 0$, $f(xt)/f(t) \to x^\alpha$ as $t \to \infty$, denoted $f \in \mathcal{R}_\alpha$.
When $1-\psi(1/\cdot) \in \mathcal{R}_{-\alpha}$ for $\alpha \in (0,1]$, it is shown in Proposition~6.1 of~\cite{Charpentier/Fougeres/Genest/Neslehova:2014} that $C_{\psi, \ell}$ is in the maximum domain of attraction of the extreme-value copula $C_{\ell_\alpha}$, i.e., for any $\bm{u} \in [0,1]^d$,
\begin{equation}\label{eq:AXCattractor}
	\lim_{n\to \infty} C_{\psi,\ell}^n(\bm{u}^{1/n}) = C_{\ell_\alpha} (\bm{u}),
\end{equation}
where for any $\bm{x} \in \mathbb{R}_+^d$, $\ell_\alpha(\bm{x}) = \ell^\alpha(\bm{x}^{1/\alpha})$. It is apparent that the Archimax family is fully flexible in the asymptotic regime, meaning that any extreme-value copula $C_\ell$ corresponds to a subclass of Archimax copulas that will be attracted to it.

\section{Model specification and notation} \label{sec:modelspecification}
	As a first step towards the specification of clustered Archimax copulas, we need to introduce the notion of clusters to the random vector $\bm{X}$. To that end, let $\mathcal{G} = \{\mathcal{G}_1,\ldots, \mathcal{G}_K\}$ be a partition of $\{1,\ldots, d\}$ into $K$ disjoint sets. Because the stochastic representation  \eqref{eq:rep} only makes sense in dimensions two and higher, we shall require, throughout this paper, that $d_k = |\mathcal{G}_k| \ge 2$ for all $k \in \{1,\ldots, K\}$.
Note that singleton clusters could be included, but this would require more tedious notation.
In this setup, $K \le \lfloor d/2 \rfloor$ and of course also $d_1 + \cdots + d_K = d$.
For convenience, we treat the subsets $\mathcal{G}_k$ as ordered sets.
This allows us to refer to the subvector of $\bm{X}$ associated with the $k$th cluster $\mathcal{G}_{k} = \{i_1,\dots,i_{d_k}\}$ as $\bm{X}_k = (X_{i_1},\dots,X_{i_{d_k}})$ or, in shorter notation, $\bm{X}_k = (X_i)_{i \in \mathcal{G}_k}$.

As we shall see shortly, a clustered Archimax copula is defined through a partition $\mathcal{G}$ as well as $K$ stdfs and $K$ distortion variables. 
To ease the reading, we will use the notation $\bm{\ell}=(\ell_1,\ldots, \ell_K)$ and $\bm{\psi}=(\psi_1,\ldots, \psi_K)$, where for each $k \in \{1,\ldots, K\}$, $\ell_k$ is a $d_k$-variate stdf and $\psi_k$ is a $d_k$-monotone Archimedean generator, i.e., differentiable up to order $d_k - 2$ with derivatives satisfying $(-1)^m\psi_k^{(m)}(x)\ge 0$ for all $x\in(0,\infty)$ for $m\in \{1,\ldots, d_k-2\}$ and further $(-1)^{d_k-2}\psi^{(d_k-2)}$ is nonincreasing and convex on $(0,\infty)$.
\begin{definition}\label{def:model}
	A $d$-variate copula $C_{\mathcal{G}, \bm{\psi}, \bm{\ell}, Q}$ is called clustered Archimax copula with cluster partition $\mathcal{G}=\{\mathcal{G}_1,\ldots, \mathcal{G}_K\}$, stdfs $\bm{\ell}$, Archimedean generators $\bm{\psi}=\{\psi_1,\ldots,\psi_K\}$ respectively $d_1,\dots,d_K$-monotone, and copula $Q$ that we term the radial copula, if it is the survival copula of a random vector $\bm{X}$ that satisfies the following:
	\begin{itemize}
		\item[(i)] For each $k \in \{1,\ldots, K\}$, $\bm{X}_k = R_k \bm{S}_k$ for some $d_k$-dimensional random vector $\bm{S}_k$ with survival function $\bar{G}_k = \bar{G}_{\ell_k}$ as in \eqref{eq:Ssurvivalfunction} and random variable $R_k$ is distributed as the inverse Williamson $d_k$-transform of $\psi_k$.
		\item[(ii)] The random vectors $\bm{S}_1,\ldots, \bm{S}_K$ are mutually independent.
		\item[(iii)] The random vector $\bm{R}=(R_1,\ldots, R_K)$ is independent of $\bm{S}_1,\ldots, \bm{S}_K$ and has copula $Q$.
	\end{itemize}
\end{definition}
	Note that by Sklar's Theorem, the joint distribution of $\bm{R}$ in Definition~\ref{def:model} is fully determined by its copula $Q$ and marginal distributions that are the inverse Williamson-$d$ transforms of $\bm{\psi}$.
	
In more explicit notation, Definition~\ref{def:model} states that, upon re-indexing so that $\mathcal{G}$ is contiguous with ordered subsets, any clustered Archimax random vector $\bm{X}$ admits the representation
\begin{equation}\label{eq:model1}
	\bm{X} = {(\bm{X}_1,\dots,\bm{X}_K) = } \bigl(R_1 S_{11},\ldots,R_1 S_{1d_1},\ldots,R_K S_{K1},\ldots,R_K S_{Kd_K}\bigr)\;.
\end{equation}
As the name suggests, certain multivariate margins of a clustered Archimax copula $C_{\mathcal{G},\bm{\psi},\bm{\ell},Q}$ are Archimax. Specifically, Theorem~3.3 of \cite{Charpentier/Fougeres/Genest/Neslehova:2014} ensures that for each $k \in \{1,\ldots, K\}$, the survival copula of $\bm{X}_k$ is the $d_k$-dimensional Archimax copula $C_{\psi_k, \ell_k}$. In particular, in the boundary case when $K=1$, the entire copula is Archimax.
It also follows from the proof of the latter theorem that the survival copula $C_{\mathcal{G},\bm{\psi},\bm{\ell},Q}$ of $\bm{X}$ in \eqref{eq:model1} is the distribution function of 
\begin{equation}\label{eq:copulamodel}
	\bigl\{\psi_1(R_1S_{11}),\ldots,\psi_1(R_1S_{1d_1}),\ldots,\psi_K(R_K S_{K1}),\ldots,\psi_K(R_K S_{Kd_K})\bigr\}\;.
\end{equation}

Throughout the paper, for each $k\in\{ 1,\ldots, K\}$, we let $\bm{S}_k=(S_i)_{i \in \mathcal{G}_k}$ be the vector that contains the components of $\bm{S}$ corresponding to cluster $k$. In particular, when $\bm{X}$ is as in \eqref{eq:model1}, then $\bm{S} = (\bm{S}_1,\dots,\bm{S}_K)$. Furthermore, analogously to $\bm{X}$ and $\bm{S}$, for any vector $\bm{x} \in \mathbb{R}_+^{d}$ and all $k \in \{1,\dots,K\}$, we denote by $\bm{x}_{k} = (x_i)_{i \in \mathcal{G}_k} \in \mathbb{R}_+^{d_k}$ the subvector of $\bm{x}$ associated with the $k$th cluster $\mathcal{G}_k$. Accordingly, for all $k \in \{1,\dots,K\}$ and $i \in \{1,\dots,d_k\}$, we use $X_{ki}$, $S_{ki}$ and $x_{ki}$ to denote the $i$th entry of $\bm{X}_k$, $\bm{S}_k$ and $\bm{x}_k$, respectively. Finally, unless otherwise stated, all operations involving one- and multi-dimensional vectors (random or not) should be understood as componentwise, e.g., for $k \in \{1,\dots,K\}$, $R_k \bm{S}_k^2 = (R_k S_i^2)_{i \in \mathcal{G}_k}$.

\section{Extremal properties}\label{sec:modelproperties}
In this section, we investigate the {extremal} behavior of a clustered Archimax copula $C_{\mathcal{G},\bm{\psi},\bm{\ell},Q}$. The main result, Theorem~\ref{thm:MaxDofA} below, delineates the conditions under which $C_{\mathcal{G},\bm{\psi},\bm{\ell},Q}$ is in a copula domain of attraction of some extreme-value copula and identifies the latter.
Since the survival copula $C_{\mathcal{G},\bm{\psi},\bm{\ell},Q}$ of $\bm{X}$ in \eqref{eq:model1} is also the copula of $1/\bm{X}$, we will study the extremal behavior of $1/\bm{X}$.

The distortion vector $\bm{R}$ has an effect on both inter- and intra-cluster dependence at extreme levels. Its extreme behavior is important, so it is natural to make the following two assumptions. The first concerns the properties of the margins of $1/\bm{R}$. Recall that a univariate random variable $X$ with distribution $F$ is in the maximum domain of attraction of a non-degenerate distribution $G$, denoted $X\in\mathcal{M}(G)$ or $F\in\mathcal{M}(G)$ iff there exist sequences $a_n\in \mathbb{R}_+$, $b_n \in \mathbb{R}$ such that, for any $x\in\mathbb{R}$, $F^n(a_nx+b_n)\to G(x)$ as $n\to\infty$. Moreover, the Fisher-Tippett Theorem states that $G$, up to location and scale, is either Fréchet ($\Phi_\rho$), Gumbel ($\Lambda$) or Weibull ($\Psi_\rho$) with $\rho>0$.
\begin{assumption}\label{assumption:MDA}
	For a clustered Archimax copula as in Definition~\ref{def:model}, assume that $\{1,\ldots,K\}$ is the {union of  disjoint sets} $\mathcal{D}_1$ and $\mathcal{D}_2$, {such that}
	\begin{itemize}
		\item[(i)] $k\in\mathcal{D}_1$ if and only if $1/R_k\in\mathcal{M}(\Phi_{\rho_k})$ for some $\rho_k\in(0,1)$.
		\item[(ii)] $k\in\mathcal{D}_2$ if and only if there exists an $\epsilon_k>0$ such that $\E \{1/R_k^{1+\epsilon_k}\}<\infty$.
	\end{itemize}
\end{assumption}
While the two cases above cover most widely considered Archimedean generators, we do conjecture an extension in Appendix~\ref{assumption:MDAExt} that includes the boundary case $1/R_k\in \mathcal{M}(\Phi_1)$ and $\E (1/R_k) = \infty$. If Assumption~\ref{assumption:MDA} holds, $k \in \mathcal{D}_1$ means that $1/R_k$ is heavy-tailed and this occurs if and only if $1-\psi(1/\cdot) \in \mathcal{R}_{-\rho_k}$, as shown in Theorem 2 in \cite{Larsson/Neslehova:2011}. In contrast, $k \in \mathcal{D}_2$ implies that $1-\psi(1/\cdot) \in \mathcal{R}_{-1}$ by Proposition~2 in \cite{Belzile/Neslehova:2017}. By the same proposition, one then has that $1/{X_{ki}} \in \mathcal{M}(\Phi_{\rho_k})$ for $k \in  \mathcal{D}_1$ and $i \in \{1,\ldots, d_k\}$ and $1/X_{ki} \in \mathcal{M}(\Phi_1)$ for $k \in  \mathcal{D}_2$ and $i \in \{1,\ldots, d_k\}$. This means that under Assumption \ref{assumption:MDA}, the respective clustered Archimax copula is in the copula domain of attraction of an extreme-value copula $C_0$ if and only if  $1/\bm{X}$ is in the maximum domain of attraction of an extreme-value distribution with copula $C_0$. Such a domain of attraction result requires further assumptions on the extremal behavior of the entire vector $1/\bm{R}$.

\begin{assumption}\label{assumption:stdf}
	For a clustered Archimax copula as in Definition~\ref{def:model}, assume that the reciprocal distortion vector $1/\bm{R}$ is in the maximum domain of attraction of a multivariate extreme-value distribution with stdf $\ell_{1/\bm{R}}$ given, for $\bm{x} \in \mathbb{R}_+^K$, by
$$
	\ell_{1/\bm{R}}(\bm{x})= \E\Bigl(\max_{1 \leq k \leq K} x_k W_k \Bigr)
$$
for some positive random variables $W_1,\ldots,W_K$ with unit mean.
\end{assumption}
It will become apparent in the next result that the choice of the aforementioned $d$-norm representations for stdfs is convenient in this context. {We are now in position to formulate the main result of this section.}
\begin{theorem}\label{thm:MaxDofA}
	Let $C_{\mathcal{G},\bm{\psi},\bm{\ell}, Q}$ be a clustered Archimax copula such that Assumptions~\ref{assumption:MDA} and \ref{assumption:stdf} hold with $(W_1,\ldots, W_K)$ independent of $\bm{S}$.
	For $k\in\mathcal{D}_1$, let {$b_k= \E\{1/Z_k^{\rho_k} \}$, $Z_k\sim \mathrm{Beta}(1,d_k-1)$}. Then $1/\bm{X} \in \mathcal{M}(H)$ with $1/X_{ki} \in \mathcal{M}(H_{ki})$, where $H_{ki}=\Phi_{\rho_k}$ for $k\in\mathcal{D}_1$ and $i \in \{1,\ldots, d_k\}$ and $H_{ki}=\Phi_{1}$ for $k\in\mathcal{D}_2$ and $i\in\{1,\ldots,d_k\}$. The  stdf of $H$ is given for all $\bm{x}\in\mathbb{R}_+^d$ by
	\begin{align}\label{eq:stdf}
	\ell_{\mathcal{G},\bm{\psi},\bm{\ell},Q} (\bm{x}) =
	\E\left( \max_{\substack{ k\in\mathcal{D}_1 \\ 1 \leq i \leq d_k}} \frac{x_{ki}W_k}{b_k S_{ki}^{\rho_k}} \right)
	+\sum\limits_{k\in\mathcal{D}_2}\ell_k(\bm{x}_{k})\;.
	\end{align}
\end{theorem}
Inter-cluster asymptotic independence can also be achieved if the distortions are asymptotically independent, as shown in the following corollary.
\begin{corollary}\label{cor:AI}
	If $\{1/R_k : k\in\mathcal{D}_1\}$ are asymptotically independent, then  for $\bm{x}\in\mathbb{R}_+^d$, the limiting stdf in \eqref{eq:stdf} simplifies to
	\begin{align*}
	\ell_{\mathcal{G},\bm{\psi},\bm{\ell},Q}(\bm{x})=&\sum\limits_{k\in\mathcal{D}_1}\ell^{\rho_k}_k\bigl(\bm{x}_k^{1/\rho_k}\bigr)+\sum\limits_{k\in\mathcal{D}_2}\ell_k(\bm{x}_k)\;.
	\end{align*}
\end{corollary}
\begin{remark}
	Note that {under the hypothesis of Theorem~\ref{thm:MaxDofA},} the asymptotic behavior of $\{1/R_k:k\in\mathcal{D}_2\}$ {has no influence on the form of }$\ell_{\mathcal{G},\bm{\psi},\bm{\ell},Q}$.
\end{remark}
The following corollary to Theorem~\ref{thm:MaxDofA} compares the inter-cluster stdf to that of the reciprocal distortions $(1/R_1,\ldots,1/R_K)$.}
\begin{corollary}\label{cor:asymptoticordering}
	Under the hypothesis of Theorem~\ref{thm:MaxDofA}, let  $(i_1,\ldots, i_K) \in \mathcal{G}_1 \times \ldots \times \mathcal{G}_K$. Then, for all $\bm{x}'\in\mathbb{R}_+^K$ and $\bm{x} \in\mathbb{R}_+^d$ such that $(x_{i_1},\dots,x_{i_K}) = \bm{x}'$ and $x_i = 0$ for all $i \in \{1,\dots,d\} \setminus \{i_1,\dots,i_K\}$,
	$$\ell_{1/\bm{R}}(\bm{x}')\le \ell_{\mathcal{G},\bm{\psi},\bm{\ell},Q}(\bm{x}).
	$$
\end{corollary}
\begin{remark}
	It is worth noting that \eqref{eq:stdf} elicits a new method to combine different stdfs in a non-trivial way. Since the second component of \eqref{eq:stdf} does not reveal any new combination of stdfs, suppose for now that $\mathcal{D}_2=\emptyset$. For a given $k\in\{1,\ldots,K\}$, we then automatically have that $\mathcal{D}_1$. Setting $x_{li}=0$ for all $l\ne k$ and all $i=1,\ldots,d_l$ recovers the marginal stdf of the cluster $k$.
This marginal stdf is equal to
$$
	\E\left(\max_{1 \leq i \leq d_k}\frac{x_{ki}}{b_k S_{ki}^{\rho_k}}\right)\;
$$	
for $\bm{x}_k \in\mathbb{R}_+^{d_k}$, which itself is equal to $\ell_k^{\rho_k}\bigl(\bm{x}_k^{1/\rho_k}\bigr)$ by Proposition~(6.1) of \cite{Charpentier/Fougeres/Genest/Neslehova:2014}. In the bivariate case, the form above is a special case of (7) in \cite{Engelke/Opitz/Wadsworth:2019}. The attractor of the bivariate Archimax copula is in particular obtained as a special case of their Proposition~1 and Equation~(6), see Sections 2.1 and 4 therein. The complete stdf, defined for all $\bm{x} \in \mathbb{R}_+^d$, by
$$
	\E\left( \max_{\substack{k\in\mathcal{D}_1 \\ 1 \leq i \leq d_k}} \frac{x_{ki}W_k}{b_kS_{ki}^{\rho_k}}\right)\;,
$$
essentially mixes the marginal cluster stdfs $\ell_1^{\rho_1}(\bm{x}_1^{1/\rho_1}),\ldots,\ell_K^{\rho_K}(\bm{x}_K^{1/\rho_K})$ with the limiting stdf of $(1/R_1,\ldots,1/R_K)$. Simply put, Corollary~\ref{cor:asymptoticordering} shows that this mixing results in a weaker asymptotic dependence between clusters than that of the reciprocal distortions $(1/R_1,\ldots,1/R_K)$, characterized by $\ell_{1/\bm{R}}$.
\end{remark}

\section{Insights into modeling}\label{sec:simulations}
In this section, we provide examples of parametric families that can be used to construct clustered Archimax copulas. Simulating from single Archimax copulas has recently garnered attention, as methods have been advanced by \cite{Mai:2022} and \cite{Ng/Hasa/Tarokh:2022}. Due to the popularity of Archimedean copulas, there is a wide array of parametric families for the distortions to choose from.
When $\psi$ is $d$-times differentiable, its inverse Williamson $d$-transform has the density, given, for $r>0$, by
$$
f_R(r)=(-1)^d\frac{r^{d-1}\psi^{(d)}(r)}{(d-1)!}\;;
$$
viz. Eq. (2) in \cite{McNeil/Neslehova:2010}.
\begin{example}[Clayton Generator]\label{example:clayton}
	Using the inverse Williamson $d$-transform, one can obtain the distribution of $R$ in the case when $\psi_\theta$ is Clayton with parameter $\theta$. In the Clayton case, as seen in \cite{McNeil/Neslehova:2009}, one has for {$r>0$},
	$$
	f_R(r) = \frac{(1+\theta r)^{-d-1/\theta}\ r^{d-1}}{(d-1)!}\prod_{j=1}^{d-1}(1+\theta j)\;.
	$$
	We can see that for any $d\ge 2$ and $\beta<d$,
	$$
	\E(1/R^\beta)=\frac{\theta^d\left\{\prod_{j=0}^{d}(1/\theta+{j}) \right\}}{(d-1)!}\int_{0}^{\infty}\frac{r^{d-1-\beta}}{(1+\theta r)^{1/\theta+d}}dr<\infty\;.
	$$
	Thus if the $k$-th cluster has a Clayton distortion, then its components are asymptotically independent from all other clusters since $k\in\mathcal{D}_2$ in Theorem \ref{thm:MaxDofA}.
\end{example}
\begin{example}[Joe generator] Recall the form of the Joe generator $\psi_\theta$ from (4.2.6) in \cite{Nelsen:2006}.
	Since $1-\psi_\theta(1/\cdot)\in\mathcal{R}_{-1/\theta}$, Theorem 2 from \cite{Larsson/Neslehova:2011} implies that $1/R\in\mathcal{M}(\Phi_{1/\theta})$.
	Thus, if the $k$-th cluster has a Joe distortion, it is asymptotically dependent with all other clusters with $j\in\mathcal{D}_1$, whose distortions $R_j$ are asymptotically dependent with $R_k$.
\end{example}

We now present synthetic examples of clustered Archimax copula based on the families presented above. By virtue of being constructed via the stochastic representation in Equation \eqref{eq:copulamodel}, random number generation from this model is straightforward. Recall that we do not require complete monotonicity of the Archimedean generators and therefore rely on the radial representation in Equation \eqref{eq:rep}. Our simulation algorithm is a simple extension of Algorithm 4.2 in \cite{Charpentier/Fougeres/Genest/Neslehova:2014}. The R code to generate the samples is provided in the supplementary materials.

\begin{algorithm} \label{alg:sim}
Let $C_{\mathcal{G}, \bm{\psi}, \bm{\ell}, Q}$ be as in Definition~\ref{def:model}. To simulate an observation $(U_1,\ldots,U_d)$ from $C_{\mathcal{G}, \bm{\psi}, \bm{\ell}, Q}$,
\begin{enumerate}
	\item Simulate a vector $(R_1,\ldots,R_K)$. This can be done by simulating a vector $(V_1,\ldots,V_K)\sim \bar{Q}$ and applying the transformations $F_{R_k}^{-1}(1-V_k)$ for each $k\in\{1,\ldots,K \}$. Following \cite{McNeil/Neslehova:2009}, for $r\in[0,\infty)$,
	$$
	F_{R_k}(r)=1-\sum_{j=0}^{d_k-2}\frac{(-1)^jr^j\psi_k^{(j)}(r)}{j!}-\frac{(-1)^{d-1}r^{d-1}\psi_{k,+}^{({d_k-1})}(r)}{({d_k-1})!}\;,
	$$
	where $\psi_{k,+}^{({d-1})}$ is the right-hand derivative of $\psi_k^{(d-1)}$. 
	\item For each $k\in\{1,\ldots,K\}$, generate an observation $\bm{S}_k=(S_{k1},\ldots,S_{kd_k})$ whose survival function is given, for any $\bm{s}\in\mathbb{R}_+^{d_k}$, by
	\begin{equation*}
		\bar{G}_{\ell_k}(\bm{s})=\Pr(S_{k1}>s_1,\ldots,S_{k d_k}>s_{d_k})=\left[\max\{0,1-\ell_k(\bm{s})\}\right]^{{d_k}-1}\;.
	\end{equation*}
	\item Construct $\bm{U}$ by setting $\bm{U}_1=\psi_1(R_1\bm{S}_1),\ldots,\bm{U}_K=\psi_K(R_K\bm{S}_K)$.
\end{enumerate}
\end{algorithm}

\begin{remark}
	In fact, the dependence structure of the distortions does not need to be defined via a copula, as long as it can be simulated from.
\end{remark}
In order to illustrate the findings of Section~\ref{sec:modelproperties}, we use Algorithm~\ref{alg:sim} to generate samples from two clustered Archimax copulas, Model A and Model B described in Table~\ref{tab:sim}; they differ only in the choice of their radial copula $Q$.
In both cases, three trivariate Archimax copulas representing three clusters are combined to form a $9$-dimensional dependence structure. The first cluster is defined by a Clayton-Gumbel Archimax copula, while the other two are defined by Joe-Gumbel Archimax copulas. The partition of $\{1,\ldots,9\}$ is $\mathcal{G}=\{\mathcal{G}_1,\mathcal{G}_2,\mathcal{G}_3 \}=\{\{1,2,3\},\{4,5,6\},\{7,8,9\} \}$. Assumption~\ref{assumption:MDA} holds with $1\in \mathcal{D}_1$ and $2,3\in \mathcal{D}_2$. Figures~\ref{fig:pairplotN} and \ref{fig:pairplotG} in Appendix~\ref{sec:additionalfigures} represent samples drawn from Model A and B.
It is clearly visible when comparing the off-diagonal 3-by-3 blocks that both models have the same intra-cluster dependence while having different inter-cluster dependence.

\begin{figure}
	\begin{subfigure}[b]{0.45\linewidth}
		\includegraphics[scale=0.5]{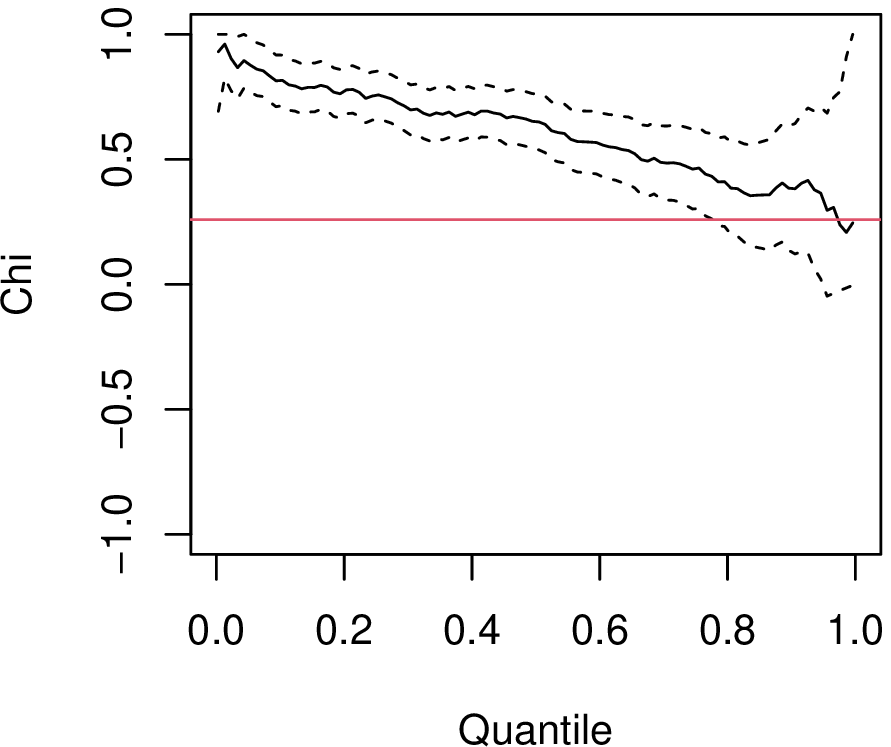}
		\caption{Model A: Pair $\{1,2\}$}
	\end{subfigure}
	\begin{subfigure}[b]{0.45\linewidth}
		\includegraphics[scale=0.5]{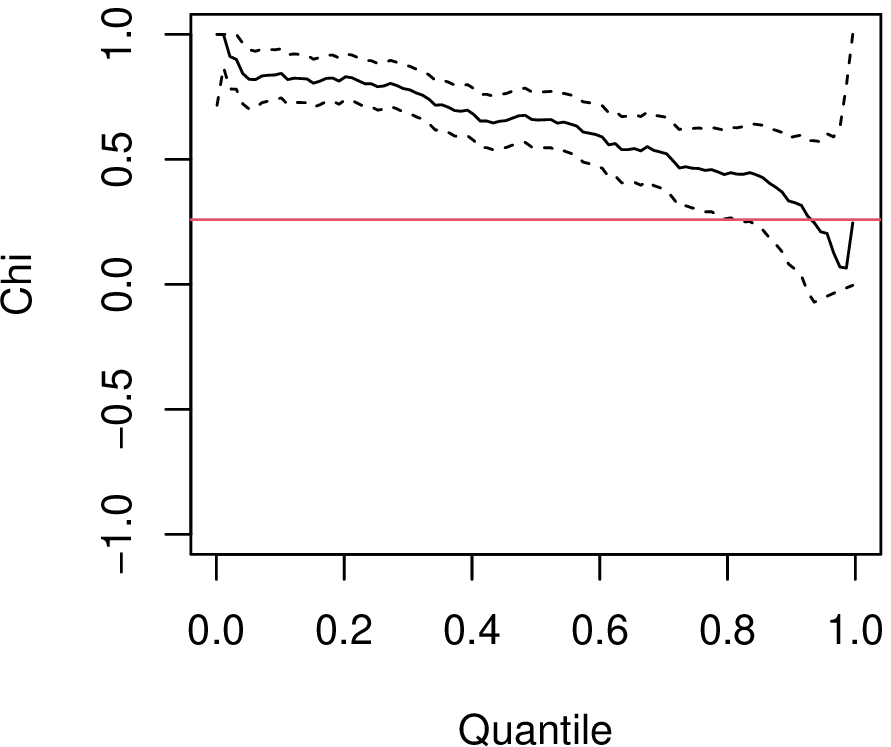}
		\caption{Model B: Pair $\{1,2\}$}
	\end{subfigure}
		\begin{subfigure}[b]{0.45\linewidth}
		\includegraphics[scale=0.5]{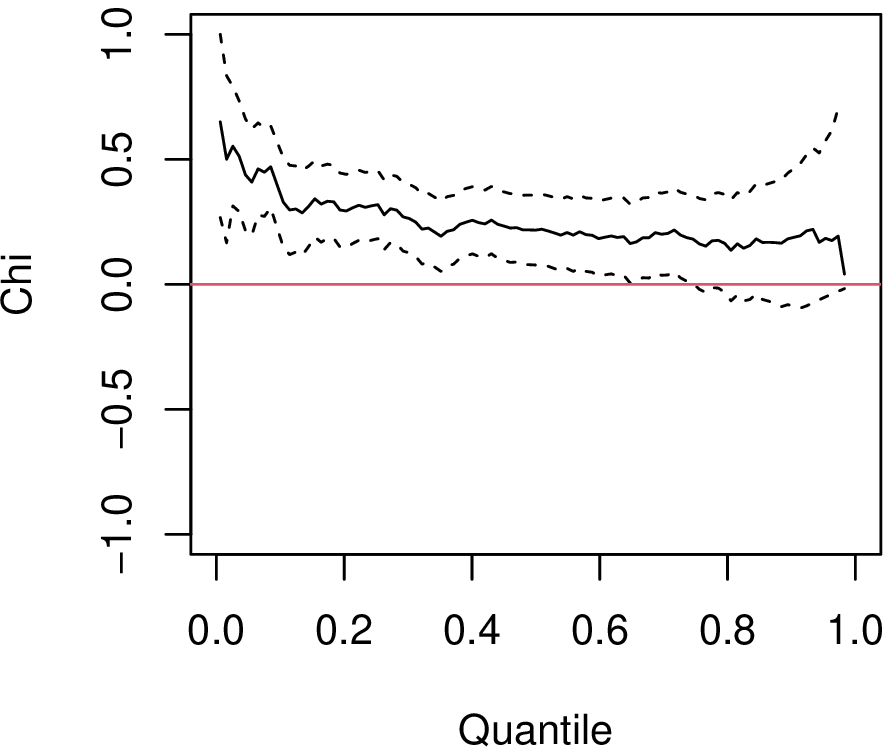}
		\caption{Model A: Pair $\{1,7\}$}
	\end{subfigure}
	\begin{subfigure}[b]{0.45\linewidth}
		\includegraphics[scale=0.5]{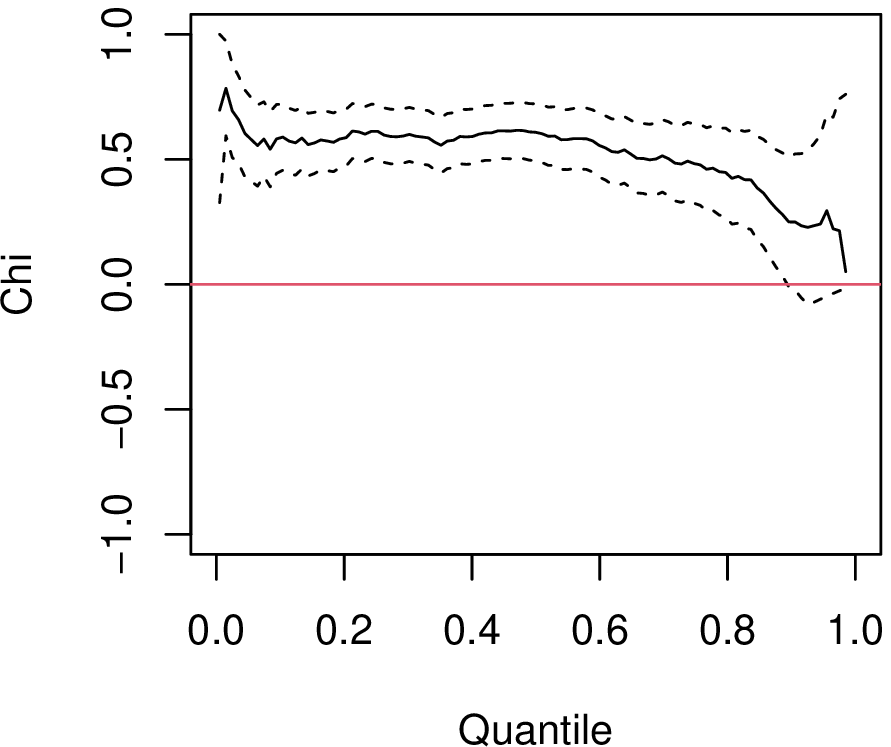}
		\caption{Model B: Pair $\{1,7\}$}
	\end{subfigure}
	\begin{subfigure}[b]{0.45\linewidth}
	\includegraphics[scale=0.5]{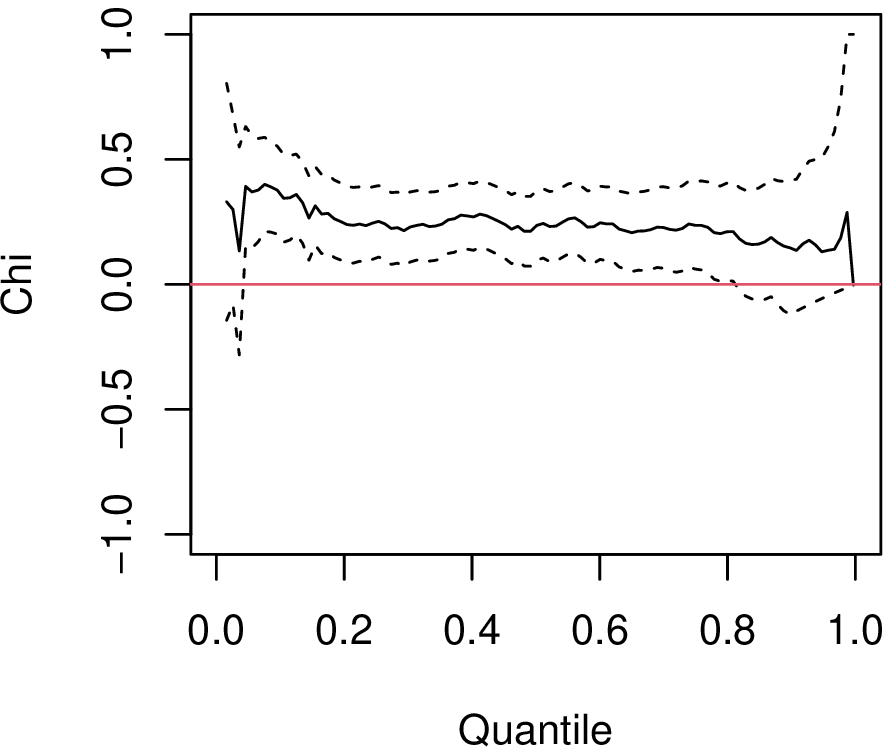}
	\caption{Model A: Pair $\{4,7\}$}
\end{subfigure}
\begin{subfigure}[b]{0.45\linewidth}
	\includegraphics[scale=0.5]{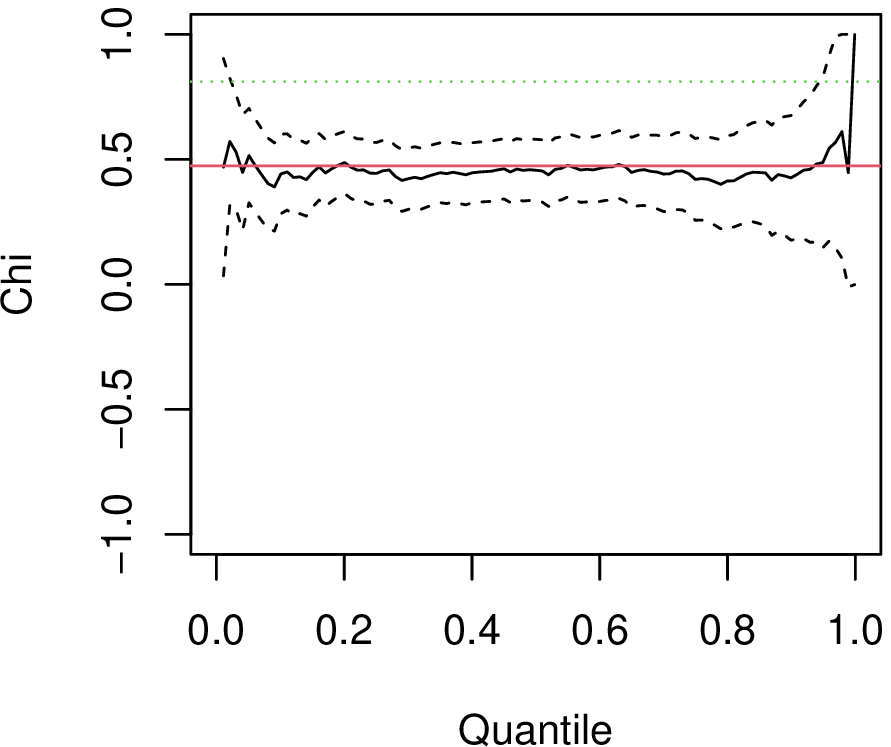}
	\caption{Model B: Pair $\{4,7\}$}
\end{subfigure}
	\caption{ Pair chi plots for the sample from Model A (left) and Model B (right). The full black line is the empirical estimate of \eqref{eq:chi}, the dotted lines are 95\% confidence intervals and the red lines represent the true values of $\lim_{q\to 1}\chi_{ij}(q)$ for each pair ${(i,j)}$. The dotted green line in panel (f) represents the upper tail dependence coefficient of $(1/R_2,1/R_3)$. The samples used for the empirical estimates are those of Figures~\ref{fig:pairplotN} and \ref{fig:pairplotG}}
	\label{fig:chiplots1}
\end{figure}

\begin{table}[ht]
	\renewcommand\arraystretch{1.5}
	\caption{Simulated models}\label{tab:sim}
	\begin{tabular}[t]{l>{\raggedright}p{0.3\linewidth}>{\raggedright\arraybackslash}p{0.3\linewidth}}
		\toprule
		&Model A&Model B\\
		\midrule
		Radial (survival) copula $\bar{Q}$&Gaussian, $\rho_{12}=\rho_{13}=\rho_{23}=0.5$& Gumbel, $\vartheta_{\bm{R}}=4$\\
		Radial extremal behavior&Asymptotic independence $\lambda_{\bm{R}}=0$&Asymptotic dependence 
		$\lambda_{\bm{R}}=0.81$\\
		$\bm{S}$ specification&
		$\ell_1$: Logistic, $\vartheta_1=1.25$ 
		
		$\ell_2$: Logistic, $\vartheta_2=2$
		
		$\ell_3$: Logistic, $\vartheta_3=1.5$&
		$\ell_1$: Logistic, $\vartheta_1=1.25$ 
		
		$\ell_2$: Logistic, $\vartheta_2=2$
		
		$\ell_3$: Logistic, $\vartheta_3=1.5$\\
		$\bm{R}$ specification&
		$R_1$: Clayton, $\theta_1=1.5$ 
		
		$R_2$: Joe, $\theta_2=1.5$
		
		$R_3$: Joe, $\theta_3=2$&
		$R_1$: Clayton, $\theta_1=1.5$ 
		
		$R_2$: Joe, $\theta_2=1.5$
		
		$R_3$: Joe, $\theta_3=2$\\
		Intra-cluster extremal dependence &
		Cluster 1: Asymptotic dependence $\lambda^{(1)}\approx 0.26$
		
		Cluster 2: Asymptotic dependence $\lambda^{(2)}\approx 0.74$ 
		
		Cluster 3: Asymptotic dependence $\lambda^{(3)}\approx 0.74$  & 
		Cluster 1: Asymptotic dependence $\lambda^{(1)}\approx 0.26$
		
		Cluster 2: Asymptotic dependence $\lambda^{(2)}\approx 0.74$ 
		
		Cluster 3: Asymptotic dependence $\lambda^{(3)}\approx 0.74$ \\
		Inter-cluster extremal dependence &
		Clusters 1-2: Asymptotic independence $\lambda^{(12)}= 0$
		
		Clusters 1-3: Asymptotic independence $\lambda^{(13)}= 0$
		
		Clusters 2-3: Asymptotic independence  $\lambda^{(23)}= 0$
		& 
		
		Clusters 1-2: Asymptotic independence $\lambda^{(12)}= 0$
		
		Clusters 1-3: Asymptotic independence $\lambda^{(13)}= 0$
		
		Cluster 2-3: Asymptotic dependence $\lambda^{(23)}\approx 0.5$ \\
		\bottomrule
	\end{tabular}
	
\end{table}
For both samples, we produce chi plots to represent extremal intra and inter-cluster dependence. Out of the 36 possible variable pairs, we chose 6 to cover all pairings of clusters $\mathcal{G}_1,\mathcal{G}_2,\mathcal{G}_3$. Figure~\ref{fig:chiplots1} displays those of variables $\{1,2\}$ (intra-cluster $\mathcal{G}_1-\mathcal{G}_1$), $\{1,7\}$ (inter-cluster $\mathcal{G}_1-\mathcal{G}_3$) and $\{4,7\}$ (inter-cluster $\mathcal{G}_2-\mathcal{G}_3$). Figure~\ref{fig:chiplots2} in Appendix~\ref{sec:additionalfigures} pertains to the remaining three cluster pairings. Following \cite{Coles/Heffernan/Tawn:1999}, for the pair $(i,j)$, the quantity of interest is given, for $q\in(0,1)$, by
\begin{equation}\label{eq:chi}
	\chi_{ij}(q)=2-\frac{\log \Pr (U_i<q,U_j<q) }{\log q}\;.
\end{equation}
The well-known upper tail dependence coefficient of \cite{Joe:2015} is then simply expressed
\begin{equation}\label{eq:lambda}
	\lambda_{ij}=\lim\limits_{q\uparrow 1}\chi_{ij}(q)\;,
\end{equation}
provided the limit exists. For each $k\in\{1,2,3\}$, the values of $\lambda_{ij}$ are in fact equal for all $i,j\in\mathcal{G}_k$ and $i\ne j$ due to the fact that the logistic stdf is symmetric with respect to permutation of its arguments. Therefore, for each $k\in\{1,2,3\}$, we can simplify the notation by having $\lambda^{(k)}= \lambda_{ij}$ with any $i,j\in\mathcal{G}_k$ such that $i\ne j$. Both samples exhibit the same intra-cluster extreme dependence, with ${\lambda}^{(1)}=2-2^{1/{\vartheta_1}}\approx 0.26$ for the first cluster (see Figures~\ref{fig:chiplots1} (a) and (b)), ${\lambda}^{(2)}=2-2^{1/{(\theta_2\vartheta_2)}}\approx 0.74$ (see Figures~\ref{fig:chiplots2} (c) and (d)) for the second cluster and ${\lambda}^{(3)}=2-2^{1/{(\theta_3\vartheta_3)}}\approx 0.74$ (see Figures~\ref{fig:chiplots2} (e) and (f)) for the third cluster. Note that the pairs $(\theta_2,\vartheta_2)$, $(\theta_3,\vartheta_3)$ were chosen to be different while resulting in the same upper tail dependence coefficient. These values are the same within each cluster because $C_{\psi_{\theta_k},\ell_{\vartheta_k}}$ was chosen to be exchangeable for each $k\in\{1,2,3\}$. The upper tail dependence coefficient of any bivariate copula $C$ in the domain of attraction of an extreme-value copula $C_\ell$ can be shown to be $\lambda=2-\ell(1,1)$.  We therefore obtain the true values of ${\lambda}^{(k)}$ via \eqref{eq:AXCattractor}, noting that the index of regular variation $\alpha$ is equal to $1$ for any Clayton generator $\psi_\theta$ and equal to $1/\theta$ for any Joe generator $\psi_\theta$.

Now, for $k,l\in\{1,2,3\}$ and $k\ne l$, let $\lambda^{(kl)}=\lambda_{ij}$ where $i\in\mathcal{G}_k$ and $j\in\mathcal{G}_l$. As for the intra-cluster extreme dependence, this simplified notation is possible due to the fact that given a pair $\mathcal{G}_k$ and $\mathcal{G}_l$, all values of $\lambda_{ij}$ such that $i\in\mathcal{G}_k$ and $j\in\mathcal{G}_l$ are equal.  In Figures~\ref{fig:chiplots1} (c) and (d) and \ref{fig:chiplots2} (a) and (b), we have $\lambda^{(12)}=\lambda^{(13)}=0$. As explained in Example~\ref{example:clayton}, the choice of a Clayton generator forces cluster 1 to be asymptotically independent from clusters 2 and 3. In Figure~\ref{fig:chiplots1} (e) we have $\lambda^{(23)}=0$. This is due to the fact that the Gaussian copula used to model the dependence between distortions forces asymptotic independence between clusters; see Example~\ref{cor:AI}. However, Figure~\ref{fig:chiplots1} (f) shows asymptotic dependence between clusters 2 and 3, i.e. $\lambda^{(23)}>0$. The value, approximately equal to $0.5$, was obtained numerically by simulating from \eqref{eq:stdf}. Corollary~\ref{cor:asymptoticordering} is illustrated by the fact that $\lambda^{(23)}$ is lower than the upper tail dependence coefficient of $(1/R_2,1/R_3)$.

\section{Inference when $\mathcal{G}$ is known} \label{sec:inference}

We first discuss inference for clustered Archimax copulas when the partition $\mathcal{G}$ is known or hypothesized. In some cases, such as a portfolio containing stocks from distinct industries, this partition can be inherent to the dataset. When $\mathcal{G}$ is unknown, one can employ clustering algorithms as discussed in Section~\ref{sec:data-application}.

From now on, let $\bm{Y}$ denote the $d$-dimensional, continuous random vector of interest, whose underlying copula $C_{\mathcal{G},\bm{\psi},\bm{\ell},Q}$ is clustered Archimax. For $k\in \{ 1,\ldots,K\}$ and $i\in\{1,\ldots,d_k\}$ let $F_{ki}$ denote the distribution function of $Y_{ki}$, recalling that $\bm{Y}_k = (Y_i)_{i \in \mathcal{G}_k}$ by definition, and let $\bm{U}$ be the $d$-dimensional copula realization associated with $\bm{Y}$, i.e., $U_{ki} = F_{ki}(Y_{ki})$ and $\bm{U}$ has distribution $C_{\mathcal{G},\bm{\psi},\bm{\ell},Q}$.
Finally, suppose that we observe a sample $\{\bm{Y}^j\}_{j=1}^n$ of i.i.d.{} replicates of $\bm{Y}$ with corresponding (unobserved) copula realizations $\{\bm{U}^j\}_{j=1}^n$.
	 
Theorem~3.3 of \cite{Charpentier/Fougeres/Genest/Neslehova:2014} implies that there exists a stochastic representation similar to \eqref{eq:model1} for the copula $C_{\mathcal{G},\bm{\psi},\bm{\ell},Q}$.
We will lean on this representation to propose a method for inferring $C_{\mathcal{G},\bm{\psi},\bm{\ell},Q}$ from $\{\bm{Y}^j\}_{j=1}^n$.
Given $\mathcal{G}$, we perform inference for $\bm{\psi}$, $\bm{\ell}$ and $Q$ separately.
To ease the reading, we denote by $\mathcal{I}$ the set
	\begin{align} \label{eq:I}
		\mathcal{I} = \{(i,j,k): 1 \leq k \leq K,\ i,j \in \mathcal{G}_k,\ i < j\}
	\end{align}
	throughout the section, as well as in Appendix~\ref{app:jackknife}. In Section~\ref{sec:inference-intra}, we briefly review how each individual cluster can be modeled using existing inference techniques for Archimax copulas. In Section~\ref{sec:app-inference-inter}, we propose a method for estimating the dependence between clusters.

\subsection{Individual cluster inference}\label{sec:inference-intra}

This section pertains to the estimation of each marginal Archimax copulas of the clustered Archimax copula $C_{\mathcal{G},\bm{\psi},\bm{\ell},Q}$. As such, inference methods have already been developed in \cite{Chatelain/Fougeres/Neslehova:2020} and more recently in \cite{Ng/Hasa/Tarokh:2022}.

\subsubsection{Estimating $\psi_1,\dots,\psi_{K}$} \label{sec:inference-psi}

We use a parametric approach to estimate the generators $\bm{\psi} = (\psi_1,\dots,\psi_{K})$, which requires selecting a (possibly different) parametric family for each generator. Certain properties observed in the data can guide the user to specific choices of generators. For example, the presence of within-cluster lower tail dependence warrants the use of a Clayton generator. However, as explained in Example~\ref{example:clayton}, this choice implies the asymptotic independence of $\bm{Y}_k$ and $\bm{Y}_\ell$ for all $k,\ell \in \{1,\dots,K\}$ with $k \neq \ell$. This consequence is desirable when studying the precipitation data in Section~\ref{sec:data-application}, but would not be appropriate, for instance, when modeling temperatures over the same region; see \cite{Davison/Huser/Thibaud:2013} for more details.

Once the parametric family of each generator is chosen, the parameters $\bm{\theta} = (\theta_1,\dots,\theta_K)$ of $\bm{\psi}$ remain to be estimated.
For a given $k \in \{1,\dots,K\}$, our estimator of $\theta_k$ relies on the following remark.
\begin{remark} \label{rem:theta-constraint}
	Any margin of an Archimax copula is itself Archimax with the same generator.
	In particular, for all $(i,j,k) \in \mathcal{I}$, with $\mathcal{I}$ as in \eqref{eq:I}, the distortions parameters $\theta_{k}$ and $\theta_{ij}$ associated with $\bm{Y}_k$ and $(Y_i,Y_j)$, respectively, are identical, i.e., $\theta_k = \theta_{ij}$.
\end{remark}
\noindent
Remark~\ref{rem:theta-constraint} suggests that for each $(i,j,k) \in \mathcal{I}$, an estimator $\hat{\theta}_{ij}$ of $\theta_{ij}$ also estimates $\theta_k$.
In view of this, we estimate $\bm{\theta}$ by $\bar{\bm{\theta}} = (\bar{\theta}_1,\dots,\bar{\theta}_K)$ given, for $k \in \{1,\dots,K\}$ and $\mathcal{G}_{kk} = \{(i,j): i,j \in \mathcal{G}_k,\ i < j\}$, by
\begin{equation} \label{eq:theta-bar}
\bar{\theta}_k = \frac{1}{\binom{d_k}{2}}\sum_{(i,j) \in \mathcal{G}_{kk}} \hat{\theta}_{ij},
\end{equation}
where, for any $(i,j) \in \mathcal{G}_{kk}$, $\hat{\theta}_{ij}$ is the estimator defined in Section~7 of \cite{Chatelain/Fougeres/Neslehova:2020}.
These latter authors discuss the cases when the underlying generator is of the Clayton, Genest-Ghoudi, Frank or Joe families.

Remark~\ref{rem:theta-constraint} also suggests a simple way of assessing whether the partition $\mathcal{G}$ is appropriate: if the underlying distribution is indeed a clustered Archimax with partition $\mathcal{G}$, then the hypothesis $H_0$ that $\theta_{ij} = \theta_{k}$ for all $(i,j,k) \in \mathcal{I}$ must necessarily hold.
We propose to test $H_0$ using the statistic $\bm{T}$ given by
\begin{equation} \label{eq:T}
	T_{ijk} = \hat{\theta}_{ij} - \bar{\theta}_k.
\end{equation}
Although its definition involves three indices, we treat $\bm{T}$ as a vector of dimension $|\mathcal{I}|$; the specific ordering of its entries does not matter, as long as it is kept fixed.
Because we expect that departures from the null will cause certain entries of $\bm{T}$ to be large in absolute values, we test $H_0$ using either the supremum norm $\| \bm{T}\|_{\infty} = \max_{\iotaa \in \mathcal{I}} |T_{\iotaa}|$ or the (squared) Euclidean norm $\| \bm{T}\|_{2}^2 = \sum_{\iotaa \in \mathcal{I}} T_{\iotaa}^2$ of $\bm{T}$; see \cite{Perreault/Neslehova/Duchesne:2022} for similar hypothesis tests, albeit in a nonparametric context.

To derive the null distribution of $\| \bm{T} \|$, where $\| \cdot \|$ is either the supremum or Euclidean norm, we exploit the fact that for any $(i,j,k) \in \mathcal{I}$, $\hat{\theta}_{ij}$ is a function of two U-statistics with square integrable kernels \citep[Section~7]{Chatelain/Fougeres/Neslehova:2020}.
Consequently, for any $\iotaa \in \mathcal{I}$, $T_{\iotaa}$ is a function of several such U-statistics.
Standard results about U-statistics \citep{Hoeffding:1948} and the delta method then imply that, under the null, $\sqrt{n}\bm{T}$ is asymptotically Normal as $n \to \infty$, i.e., $\sqrt{n}\bm{T} \rightsquigarrow \mathcal{N}(\bm{0},\bm{\Sigma})$ for some positive definite matrix $\bm{\Sigma}$ of appropriate dimensions.
In view of Slutsky's Lemma and the Continuous Mapping Theorem, a test of approximate level $\alpha \in (0,1)$ then consists in rejecting $H_0$ whenever $\Pr(\|\bm{Z}\|  > \|\bm{T}\|) < \alpha$, where $\bm{Z} \sim \mathcal{N}(\bm{0},\hat{\bm{\Sigma}})$ for some consistent estimator $\hat{\bm{\Sigma}}$ of $\bm{\Sigma}$.
In the data application of Section~\ref{sec:data-application}, we use a jackknife estimator of $\bm{\Sigma}$ derived from Section~2(c) of \cite{Arvesen:1969}; see Appendix~\ref{app:jackknife} for the implementation details.

\subsubsection{Estimating $\ell_1,\dots,\ell_K$} \label{sec:inference-S}

To estimate the stdfs $\bm{\ell} = (\ell_1,\dots,\ell_K)$, we exploit the fact that for any $k \in \{1,\dots,K\}$ and $\bm{x} \in \mathbb{R}^{d_k}$, $\ell_k(\bm{x}) = \| \bm{x} \| A_k(\bm{x}/\|\bm{x}\|)$, where $A_k$ is the Pickands dependence function associated to $\bm{Y}_k$. Specifically, we replace $A_k$ in this latter equation by its CFG-type estimator \citep{Caperaa/Fougeres/Genest:1997}, defined as follows.
For $j \in \{1,\dots,n\}$, let $\hat{\bm{U}}^j$ be such that, for all $k \in \{1,\dots,K\}$ and $i \in \{1,\dots,d_k\}$,
\begin{equation} \label{eq:Uhat}
	\hat{U}_{ki}^j = r_{kij}/(n+1)\;,
\end{equation}
where $r_{kij}$ is the rank of $Y_{ki}^{j}$ among $Y_{ki}^{1},\dots,Y_{ki}^{n}$.
Now, for $\bm{w} \in \Delta_{d_k}$ with $\Delta_{d_k}$ as in Definition~\ref{def:axc}, let
$$
\hat{\xi}_{k,j}(\bm{w}) = \min_{1 \leq i \leq d_k} \phi_{\bar{\theta}_k}(\hat{U}_{ki}^{j})/w_{i}.
$$
Then, for all $k \in \{1,\dots,K\}$ and $\bm{w} \in \Delta_{d_k}$, the (endpoint-corrected) CFG-type estimator of $A_{k}$ is given by
\begin{equation} \label{eq:A-CFG}
	\log \hat{A}_k^{\rm CFG}(\bm{w}) = \frac{1}{n}\sum_{j=1}^n \left\{\log \phi_{\bar{\theta}_k}\left(\frac{j}{n+1}\right) - \log \hat{\xi}_{k,j}(\bm{w}) \right\}.
\end{equation}
We refer the reader to Section~3 of \cite{Chatelain/Fougeres/Neslehova:2020} for more details about the CFG-type estimator, as well as an alternative estimator based on that of \cite{Pickands:1981}. The asymptotic behavior of the CFG and Pickands-type estimators is established under regularity conditions on both $\psi$ and $A$; see Section~4 in \cite{Chatelain/Fougeres/Neslehova:2020}.

\subsection{Inference for the dependence between distortions}\label{sec:inference-CR}

To model the dependence between the components of $\bm{R}$, we suggest a parametric approach on the underlying copula $Q$. To this end, we make the assumption that $Q$ belongs to a parametric family $ \{ C_{\bm\vartheta}: \bm\vartheta \in \bm\varTheta \}$ of $K$-dimensional copulas, where the parameter space $\bm\varTheta$ is of arbitrary dimension. We further assume the existence of a multivariate density for $\bm{R}$ in order to proceed with a likelihood-based method. While we cannot observe $\bm{R}$ directly, we may still derive the corresponding likelihood based on pseudo-observations.

We follow \eqref{eq:copulamodel} and define the $d$-dimensional random vector $\bm{V}$, for all $k \in \{1,\dots,K\}$ and $i \in \{1,\dots,d_k\}$, by $V_{ki} = \psi_k(R_k S_{ki})$, so that $\bm{V}$ has distribution $C_{\mathcal{G},\bm{\psi},\bm{\ell},Q}$.
Now, fix $\bm{\iota} = (i_1,\ldots,i_K) \in \mathcal{G}_{\bullet} = \mathcal{G}_1 \times \ldots \times \mathcal{G}_K$ and note that since $X_{\iota_k} = R_k S_{\iota_k}$ for each $k \in \{1,\dots,K\}$, the density of the subvector $\bm{V}_{\bm{\iota}} = \{\psi_1(X_{\iota_1}),\ldots,\psi_K(X_{\iota_K})\}$ of $\bm{V}$ can be expressed, for $\bm{v} \in (0,1)^K$, as
\begin{equation} \label{eq:fV}
	f_{\bm{V}_{\bm{\iota}}}(\bm{v}) = \int_{(0,1)^K} f_{\bm{R}}\{\bm{\phi}(\bm{v})/\bm{s}\} \left\{ \prod_{k=1}^K \frac{f_{S_k}(s_k) \phi_k'(v_k)}{s_k} \right\} \mathrm{d} \bm{s}
\end{equation}
where $\bm{\phi}(\bm{v})/\bm{s} = \{ \phi_1(v_1)/s_1,\ldots,\phi_K(v_K)/s_K \}$ and, for $k \in\{1,\ldots,K\}$, $f_{S_k} $ is a $ \mathrm{Beta}(1,d_k-1)$ density.
Also note that, due to Sklar's Theorem, the joint density of $\bm{R}$ can be written, for $\bm{r} \in \mathbb{R}_+^K$ and $q$ the copula density of $\bm{R}$, as
\begin{equation} \label{eq:fR}
f_{\bm{R}}(\bm{r})=q\{F_{R_1}(r_1),\ldots, F_{R_K}(r_k)\}\left\{ \prod_{k=1}^K f_{R_k}(r_k) \right\}\;.
\end{equation}

Now, let $\mathcal{Y}_n = (\bm{Y}^j)_{j=1}^n$ and recall that $\bm{U}^j$, $j\in\{1,\ldots,n\}$, and $\bm{V}$ are equal in distribution.
The pseudo-copula observations can be plugged into the (pseudo) composite marginal log-likelihood \citep{Varin/Reid/Firth:2011} associated with the density of $\bm{V}_{\bm{\iota}}$ and given by
\begin{equation}\label{eq:loglikelihood}
	\mathcal{L}_K(\bm\vartheta|\mathcal{Y}_n)=\sum_{j=1}^{n} \sum_{\bm{i} \in \mathcal{G}_{\bullet}} \log f_{\bm{V}_{\bm{\iota}}}(\hat{\bm{U}}_{\bm{i}}^{j})\;,
\end{equation}
where $\hat{\bm{U}}_{\bm{i}}^j = (\hat{U}_{i_1}^j,\ldots,\hat{U}_{i_K}^j)$ with $\hat{\bm{U}}^j$ as in \eqref{eq:Uhat}, and $f_{\bm{V}_{\bm{\iota}}}$ as in \eqref{eq:fV}.

Due to the multi-dimensional integral in the expression of $f_{\bm{V}_{\bm{\iota}}}$ and the possibly large cardinality of $\mathcal{G}_{\bullet}$, \eqref{eq:loglikelihood} will often be difficult to compute in practice.
When the chosen parametric family for $Q$ allows it, a natural approach for reducing the computational burden is to consider lower-dimensional margins to form the composite likelihood function.
For example, when $Q$ is a Gaussian copula, it seems reasonable to use the pairwise marginal likelihood 
\begin{equation} \label{eq:loglikelihood-2}
\mathcal{L}_2(\bm{\vartheta}|\mathcal{Y}_n) = \sum_{j=1}^n \sum_{1 \leq k < \ell \leq K} \sum_{\bm{i} \in \mathcal{G}_{k\ell}} \log f_{V_{\iota_k}, V_{\iota_\ell}}(\hat{\bm{U}}_{\bm{i}}^j)\;,
\end{equation}
where $\mathcal{G}_{k\ell} = \mathcal{G}_k \times \mathcal{G}_\ell$, $\hat{\bm{U}}_{\bm{i}}^j = (\hat{U}_{i_1}^j,\hat{U}_{i_2}^j)$, and $f_{V_{\iota_k},V_{\iota_\ell}}$ is the corresponding bivariate analogue of $f_{\bm{V}_{\bm{\iota}}}$.
For more details about composite likelihood estimation, we refer the reader to \cite{Cox/Reid:2004} and \cite{Varin/Reid/Firth:2011}.

For some specific models, the similar yet simpler strategy underlying the estimation of $\bm{\theta}$ in \eqref{eq:theta-bar} can also be used.
For example, let $Q$ be a Gaussian copula with correlation matrix $(\vartheta_{k\ell})_{1 \leq k,\ell \leq K}$, so that $\bm{\vartheta}$ is of length $\binom{K}{2}$. 
For each $k,\ell \in \{1,\dots,K\}$ with $k < \ell$, one can first independently consider the bivariate log-likelihoods associated with $\mathcal{G}_{k\ell}$ in \eqref{eq:loglikelihood-2} to compute a set of pairwise parameters $\{\hat{\rho}_{\bm{i}}: \bm{i} \in \mathcal{G}_{k\ell}\}$.
These can then be averaged over to obtain a final estimate $\bar{\rho}_{k\ell}$ for $\vartheta_{k\ell}$.
This is the approach we opt for in the application of Section~\ref{sec:data-application}.

\begin{remark}
	The composite log-likelihoods $\mathcal{L}_K(\bm\vartheta|\mathcal{Y}_n)$ and $\mathcal{L}_2(\bm\vartheta|\mathcal{Y}_n)$ in \eqref{eq:loglikelihood} and \eqref{eq:loglikelihood-2}, respectively, depend on $\bm\vartheta$ through $q$ in \eqref{eq:fR} (or its bivariate analogues in the case of $\mathcal{L}_2(\bm\vartheta|\mathcal{Y}_n)$).
	The functions $f_{R_k}$, $F_{R_k}$ and $\phi_k$, $k\in\{1,\ldots, K\}$, are assumed to be known.
	In practice, one can plug in the parametric estimates proposed in Section~\ref{sec:inference-psi}.
	While outside of the scope of this work, one could consider the feasibility of estimating the distortion parameters together with $\bm\vartheta$ via pseudo-maximum composite likelihood.
\end{remark}

\section{Data illustration} \label{sec:data-application}

\begin{figure}[t]
\centering
\begin{minipage}{.4\linewidth}
\includegraphics[width=1\linewidth]{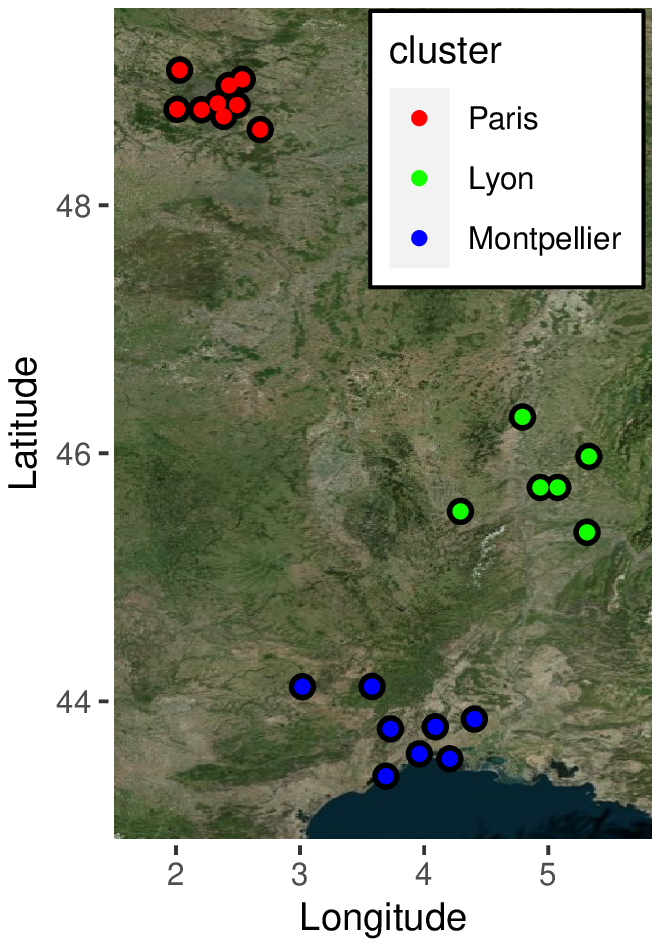}
\end{minipage}
\hspace{.05\linewidth}
\begin{minipage}{.4\linewidth}
\includegraphics[width=1\linewidth]{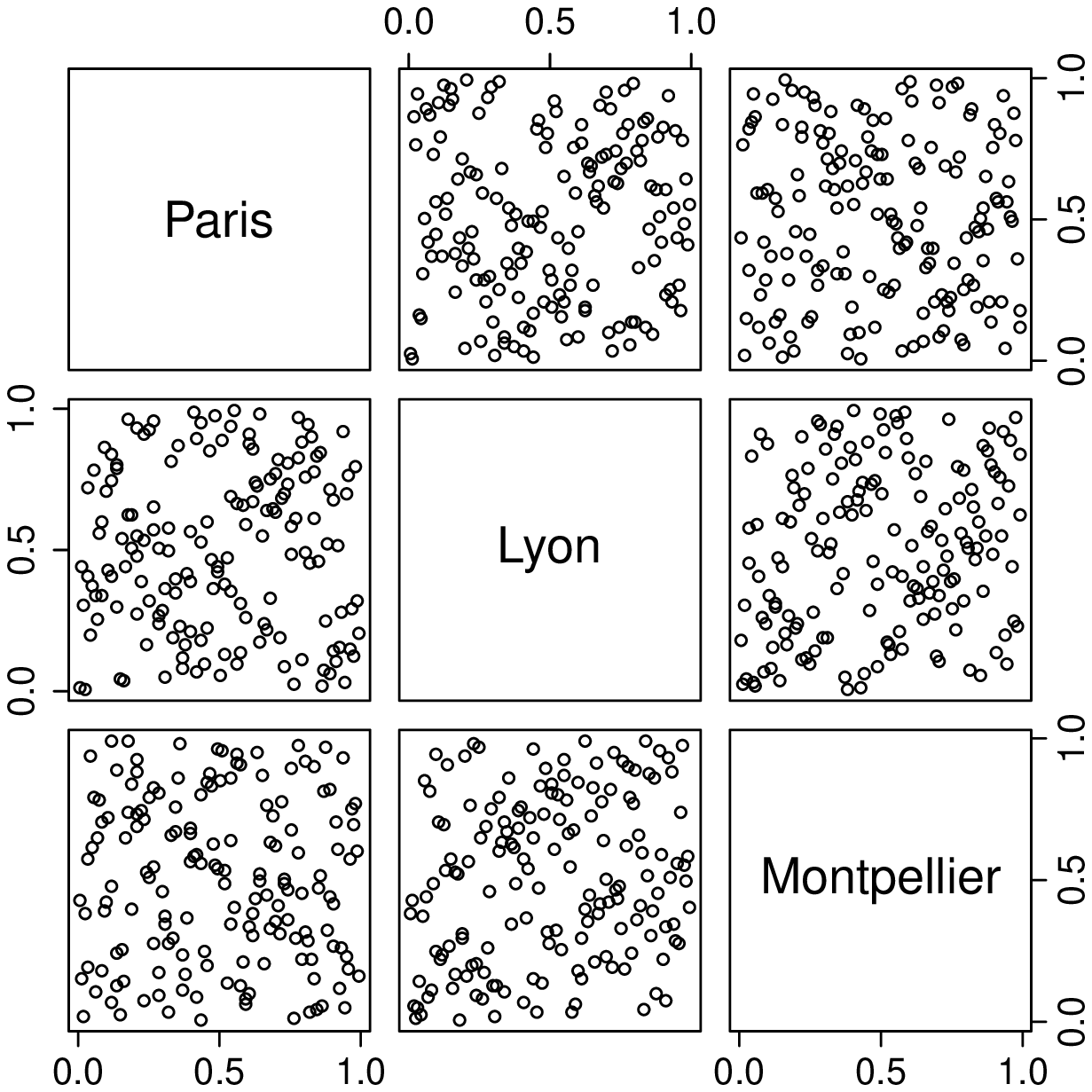}
\end{minipage}
\caption{Left: Geographical position of the stations used in the application of Section~\ref{sec:data-application}, along with their respective cluster label. Right: Pair plots of the scaled componentwise ranks of monthly maxima of precipitations for the stations around which the clusters were formed.} \label{fig:map}
\end{figure}

In this section, we illustrate the proposed methodology through an application to flood monitoring.
The data, provided by M\'{e}t\'{e}o France, consists of daily precipitation amounts measured from 1976 to 2015, inclusively, at $d=23$ meteorological stations in France.
As shown in Figure~\ref{fig:map}, the stations are agglomerated into three clusters centered around the cities of Paris ($\mathcal{G}_1$, 9 stations), Lyon ($\mathcal{G}_2$, 6 stations) and Montpellier ($\mathcal{G}_3$, 8 stations).
We thus let $\mathcal{G} = \{\mathcal{G}_1,\mathcal{G}_2,\mathcal{G}_3\}$ be the partition underlying our model.
When no such partition stems from the context of the application naturally, one can use the clustering techniques of, e.g., \cite{Bernard/Naveau/Vrac/Mestre:2013} or \cite{Saunders/Stephenson/Karoly:2021} to create a set of candidate partitions of distinct sizes.
The number of clusters can then be settled heuristically with the help of standard clustering tools (e.g., a dendrogram), or by selecting the coarsest partition that yields a satisfactory fit.

Note that our application is purposely similar to that of \cite{Chatelain/Fougeres/Neslehova:2020}, who fitted an Archimax copula to precipitation amounts (monthly maxima) measured at three nearby stations.
The greater flexibility of clustered Archimax copulas allows us to increase the number of variables considered.

\subsection{Data preprocessing} \label{sec:data}

A preliminary analysis of the data reveals the presence of seasonality and temporal dependence within the univariate series.
To mitigate the effect of seasonality, we consider only the observations from the months of September to December, inclusively, which encompass most of the extreme precipitation events.
Since our primary focus is on extreme precipitations, we then take monthly maxima of the series, yielding a total of $n=160$ observations per station.
The resulting series show no obvious sign of temporal dependence according to the Ljung-Box tests \citep{Ljung/Box:1978}; note however that the asymptotic results of \cite{Chatelain/Fougeres/Neslehova:2020} upon which the present work relies hold for alpha-mixing sequences, meaning that temporal dependence vanishing with increasing lag is indeed allowed.

The pair plots of the scaled componentwise ranks of monthly maxima of precipitations involving only the most central station of each cluster are displayed in the right panel of Figure~\ref{fig:map}.
They suggest a very weak dependence between precipitation maxima in Lyon and Montpellier, and an even weaker one, if at all, between any of these and the stations in Paris.
In contrast, similar plots for all pairs of stations within the same cluster (Figures~\ref{fig:paris}--\ref{fig:montpellier}) indicate much stronger dependencies.
In these latter plots, we also note the presence of asymmetry; this is particularly pronounced in Figure~\ref{fig:paris} (Paris).

As a final preliminary step, we used the procedure of \cite{Kojadinovic/Segers/Yan:2011} to test the hypothesis that the copula underlying each cluster of variables is an extreme-value copula.
In all three cases, the test clearly rejects the hypothesis.
This may be explained by the presence of masses of points near the bottom-left corner of many pair plots, combined with the fact that extreme-value distributions cannot allow lower-tail dependence.
In contrast with extreme-value distributions, clustered Archimax copulas may indeed allow for both lower-tail and extremal dependence; in particular, letting $\psi_{\theta}$ be the Clayton generator leads to (pairwise) lower-tail dependence coefficients equal to ${2A(1/2,1/2)}^{-1/\theta}$, where $A$ is the Pickands dependence function characterizing the Archimax, as explained in \cite{Chatelain/Fougeres/Neslehova:2020}.

\subsection{Inference for $\bm{\psi}$ and $\bm{\ell}$} \label{sec:app-inference}

The apparent lack of inter-cluster dependence in the right panel of Figure~\ref{fig:map} and the presence of lower-tail dependence in the data suggest that the Clayton generators could be good candidates for modeling the distortions, as these would produce asymptotically independent clusters.
We thus begin by formally testing whether $\mathcal{G}$ defines three asymptotically independent clusters.
To do so, we apply the test of independence for random vectors proposed by \cite{Kojadinovic/Holmes:2009}.
Because we are interested in asymptotic independence, we apply the test procedure not to our dataset of monthly maxima, but to the corresponding dataset of yearly maxima, which includes 40 observations per station;
the corresponding matrix of empirical Kendall's $\tau$ correlations is depicted in the left panel of Figure~\ref{fig:rho-tau-hat} in the Appendix.
The test yields a p-value of 0.16 for the global hypothesis of independence between the three clusters and a p-value above 0.35 for each of the three hypotheses of pairwise independence.
Although the fact that only $40$ observations are available might arguably yield a test with limited power, we move on with our analysis assuming that the clusters are asymptotically independent and that Clayton generators are reasonable choices.

\begin{figure}[t]
\centering
\includegraphics[width=.48\linewidth]{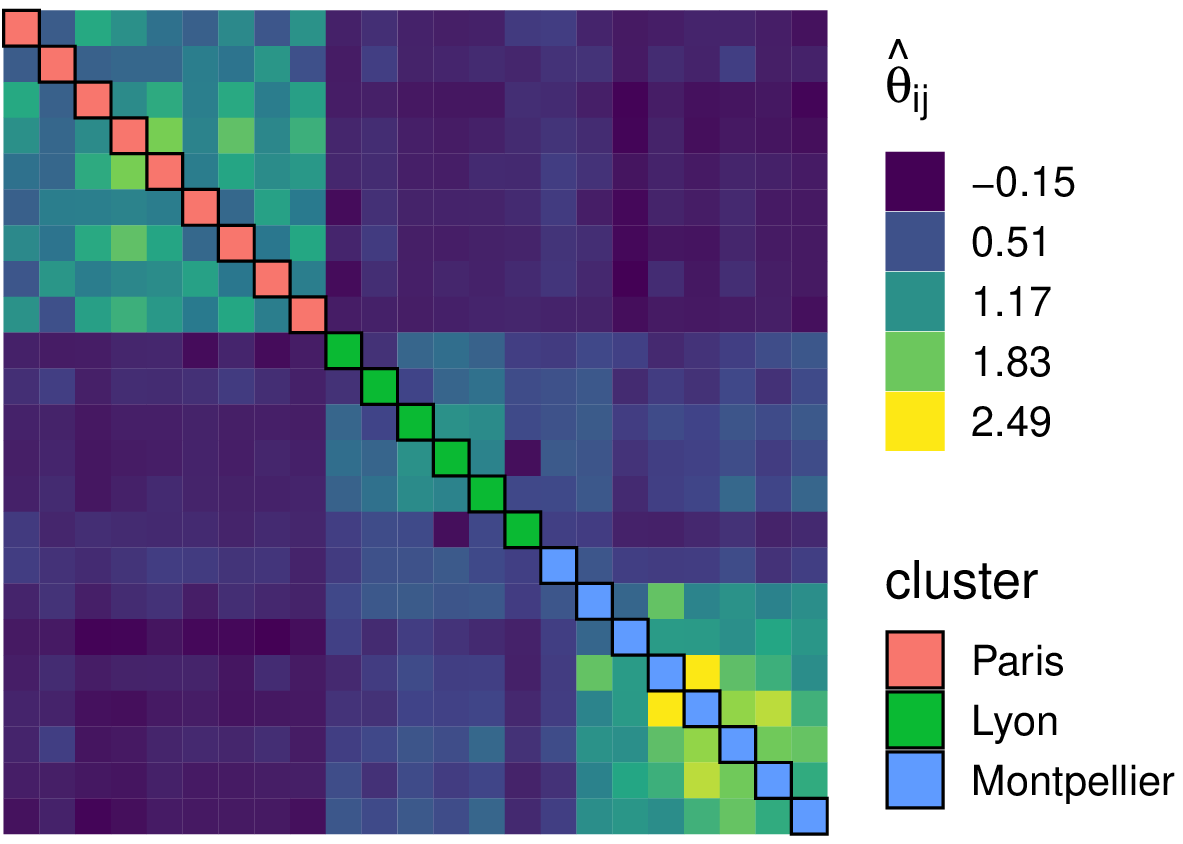}
\quad
\includegraphics[width=.48\linewidth]{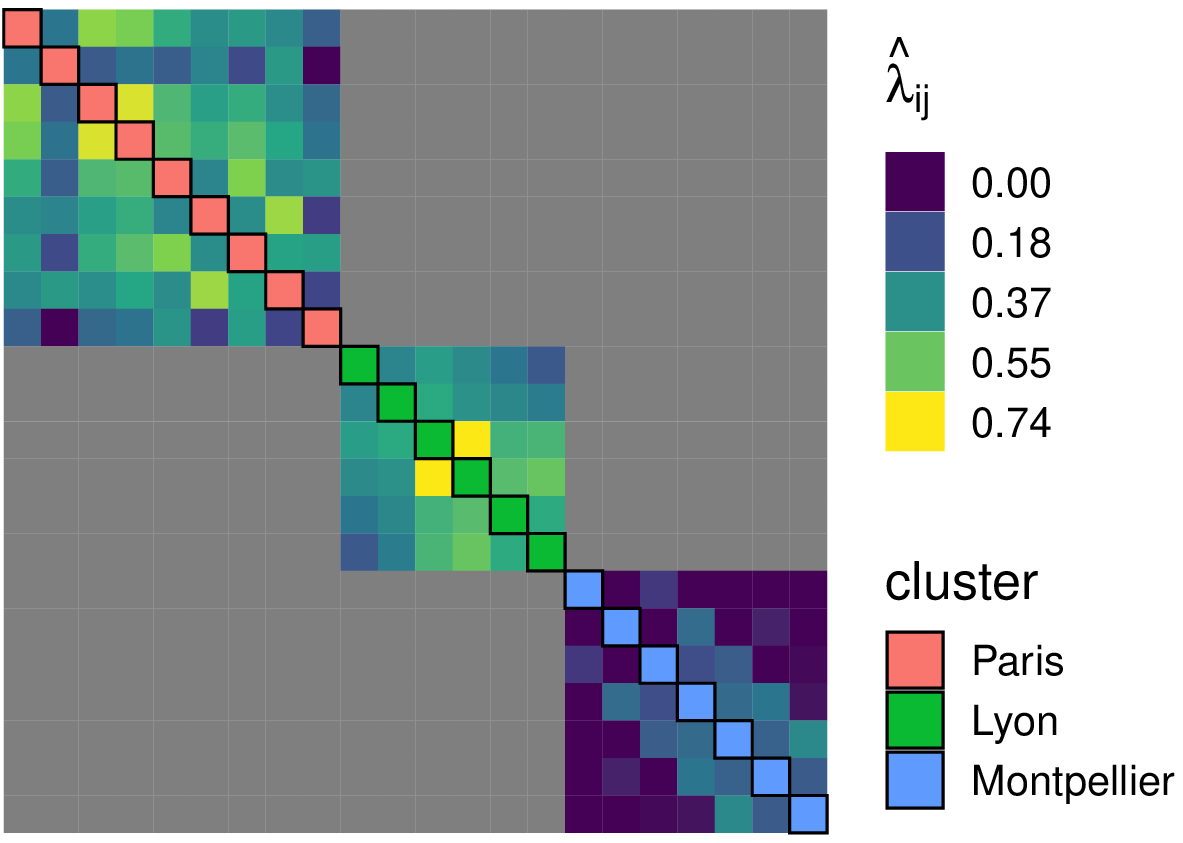}
\caption{Estimated quantities from Section~\ref{sec:app-inference}. Left: Matrix of intra-cluster pairwise upper tail dependence coefficients ($d=22$ stations). Right: Matrix of pairwise distortion parameters ($d=23$ stations).} \label{fig:theta-lambda-hat}
\end{figure}

The next step is to estimate the parameters $\theta_1$, $\theta_2$ and $\theta_3$ of the three Clayton generators.
This involves computing, for each $k \in \{1,2,3\}$, the $\binom{d_k}{2}$ pairwise estimates $\{\hat{\theta}_{ij}: i, j \in \mathcal{G}_k,\ i < j\}$; these are gathered in a $d_k \times d_k$ matrix illustrated in the left panel of Figure~\ref{fig:theta-lambda-hat}.
The resulting distortion estimates, defined in \eqref{eq:theta-bar}, are $\bar{\theta}_1 \approx 1.08$ (Paris), $\bar{\theta}_2\approx 0.62$ (Lyon) and $\bar{\theta}_3\approx 1.20$ (Montpellier), suggesting that the hypothesis of different distortions affecting the three clusters is reasonable.

At this point, one can already suspect a problem with the Montpellier cluster, as there seems to be strong discrepancies among the entries of its corresponding matrix in the left panel of Figure~\ref{fig:theta-lambda-hat}, violating the statement in Remark~\ref{rem:theta-constraint}.
To check this more formally, we perform the test described in Section~\ref{sec:inference-psi} for the hypothesis $H_0$ that $\theta_{ij} = \theta_{k}$ for all $k \in \{1,2,3\}$ and $i,j \in \mathcal{G}_k$ such that $i < j$.
While the version of the test based on the supremum norm yields an acceptable p-value (approximately $0.18$), the version based on the Euclidean norm yields an approximate p-value of $0.03$, thus rejecting $H_0$ at nominal level $0.05$.
Three similar cluster-specific tests, each involving only the pairwise distortion estimates from a single cluster, indeed reveal an anomaly with Montpellier; its corresponding approximate p-values are $0.06$ (supremum norm) and $<0.01$ (Euclidean norm).
A last series of entry-specific tests, each involving a single pairwise distortion estimate from the Montpellier cluster, strongly suggests that one particular station, the west-most station of the cluster, is at the root of the rejection.
Concretely, the $d_3-1=7$ p-values associated with this latter station are below the nominal level 0.05, with most of them being smaller than $10^{-3}$.
From the left panel of Figure~\ref{fig:map}, one can see that it effectively appears isolated from the other stations of the Montpellier cluster, as it is located on the other side of the C\'{e}vennes montain range.
Removing its corresponding column from the data yields a new distortion parameter $\bar{\theta}_3 \approx 1.48$ for the Montpellier cluster, as well as satisfactory p-values for the tests of $H_0$ and the preceeding test of asymptotic independence.
We thus remove the problematic station from the data and redefine $\mathcal{G}$ to be the ensuing partition of $\{1,\dots,22\}$, assuming that the variables were re-indexed.
Similarly, we now use $\bar{\theta}_3$ to refer to the estimate of $\theta_3$ computed without the problematic station.

Given the new partition $\mathcal{G}$, we then construct, for each $k \in \{1,2,3\}$, the semi-parametric estimate of the stdf $\ell_k$ based on $\hat{A}_k^{\rm CFG}$ defined in \eqref{eq:A-CFG}. Note that for the estimated values of $\theta_1,\theta_2,\theta_3$, Conditions 4.1 and 4.2 in \cite{Chatelain/Fougeres/Neslehova:2020} used for convergence of $\hat{A}_k^{\rm CFG}$ are met.
To better visualize the strength of the dependence between stations of the same cluster, we also computed, for each $k \in \{1,2,3\}$ and all $i,j \in \mathcal{G}_k$ such that $i \neq j$, pairwise estimates $\hat{\lambda}_{ij}$ of the upper tail dependence coefficients $\lambda_{ij}$ defined in \eqref{eq:lambda}.
To this end, we note that when Clayton generators are used, $\alpha = 1$ in \eqref{eq:AXCattractor} and the estimator in \eqref{eq:A-CFG} can be slightly modified to obtain an estimator of $\lambda_{ij}$ as follows.
Recall the definition of $\mathcal{I}$ in \eqref{eq:I}.
For all $(i,j,k) \in \mathcal{I}$, let $\ell_{ij}^*(1,1)=\ell_k(\bm{e}^{(ij)})$, where $\bm{e}^{(ij)} \in \mathbb{R}^{d_k}$ has only two non-zero entries, those associated with the $i$th and $j$th stations, equal to one.
Then, for any $(i,j,k) \in \mathcal{I}$, one can estimate $\lambda_{ij}$ combining the identities $\lambda_{ij} = 2-\ell_{ij}^*(1,1)$ and $\ell_k(\bm{x}) = \| \bm{x} \| A_k(\bm{x}/\|\bm{x}\|)$.
Because the bivariate margins of an Archimax copula are also Archimax (Remark~\ref{rem:theta-constraint}), we get that $\lambda_{ij} = 2 - 2 A_{ij}(1/2,1/2)$ with $A_{ij}(1/2,1/2) = A_k(\bm{e}^{(ij)}/2)$.
This leads to an estimator of ${\lambda}_{ij}$ given by $\hat{\lambda}_{ij} = 2 - 2 \hat{A}_{ij}^{\rm CFG}(1/2,1/2)$.
The resulting estimates are shown in the right panel of Figure~\ref{fig:theta-lambda-hat}.
Although this is not required, we expect many of these estimates to be large; this is clearly the case for the Paris and Lyon clusters, and still true, perhaps to a lesser extent, for the Montpellier cluster.

\subsection{Inference for $Q$} \label{sec:app-inference-inter}

To capture the dependence between the components of $\bm{R} = (R_1,R_2,R_3)$, we suppose that their copula $Q$ belongs to the family of Gaussian copulas with correlation parameters $\bm{\vartheta} = (\vartheta_{12}, \vartheta_{13}, \vartheta_{23})$, which we estimate by $(\bar{\rho}_{12},\bar{\rho}_{13},\bar{\rho}_{23})$, as defined at the end of Section~\ref{sec:inference-CR}. The resulting estimates are given by $\bar{\rho}_{12} \approx 0.14$ for the Paris-Lyon pair, $\bar{\rho}_{13} \approx -0.02$ for the Paris-Montpellier pair and $\bar{\rho}_{23} \approx 0.45$ for Lyon-Montpellier pair; their corresponding pairwise estimates, denoted $\{\hat{\rho}_{\bm{i}}:\bm{i} \in \mathcal{G}_k \times \mathcal{G}_\ell\}$ in Section~\ref{sec:inference-CR}, are shown in the right panel of Figure~\ref{fig:rho-tau-hat} in the Appendix.
The results suggest in particular that the clusters centered around Paris and Montpellier are nearly independent, which is coherent with the fact they are the most distant pair of clusters among the three.

This application to rainfall illustrates some of the strengths of clustered Archimax copulas. Namely, we are able to model the dependence of multi-dimensional data that consists of groups of variables which are asymptotically dependent. Within each cluster, dependence is modeled in a fully flexible way, allowing for asymmetry. Inter-cluster dependence is modeled more parsimoniously while maintaining the possibility for both asymptotic dependence and independence. Finally, the proposed copula family allows to fit data at a pre-asymptotic level while boasting flexibility at the extreme level.

\section{Discussion} \label{sec:discussion}
The clustered Archimax model studied in this paper is related to several other articles in the literature. Hierarchical constructions based on Archimax copulas were proposed by \cite{Hofert/Huser/Prasad:2017}. Specifically, their construction is based on the frailty representation of Archimax copulas, which only holds for completely monotone generators. Hierarchies can be induced via the frailties, the stdf, or both. It would be interesting to establish the attractor of their proposed hierarchical Archimax copula and compare it to that of the clustered Archimax copula. The extremal dependence structure of Liouville copulas is established in \cite{Belzile/Neslehova:2017}. The stochastic representation of Liouville copulas is similar to that of Archimax copulas, as they are survival copulas of vectors of the form $R\bm{D}$, with $R$ a nonnegative random variable and $\bm{D}$ a Dirichlet random vector. The work presented in this paper differs from this by replacing the Dirichlet component by a vector $\bm{S}$ characterized by an stdf and by allowing for multiple distorting random variables $R_1,\ldots,R_K$, thus inducing a hierarchy (or clustering). Finally, \cite{Engelke/Opitz/Wadsworth:2019} establish the extremal dependence of bivariate vectors of the form $R\times(W_1,W_2)$ for an extensive combination of asymptotic behaviors of both $R$ and $(W_1,W_2)$. The attractor of the bivariate Archimax copula is in particular obtained as a special case of their Proposition 1 and Equation (6), see Sections 2.1 and 4 therein.

\bibliographystyle{imsart-nameyear}
\bibliography{bibliography}

\numberwithin{figure}{section}
\renewcommand\thefigure{\thesection.\arabic{figure}} 
\begin{appendices}

\section{Proofs of Section~\ref{sec:modelproperties}}\label{sec:proofsofmodelproperties}
This section contains the proofs of the results {from} Section~\ref{sec:modelproperties}. {We begin with auxiliary results in Section \ref{sec:auxiliaryresults}; Theorem~\ref{thm:MaxDofA} and its Corollaries are proved in Sections \ref{sec:proofmda} and \ref{sec:sec:proofcor}, respectively.} 

\subsection{Auxiliary results}\label{sec:auxiliaryresults}
The following proposition is used to prove Theorem \ref{thm:MaxDofA} but is also of independent interest.
\begin{proposition}\label{prop:stdfS}
	Let $\bm{S}=(S_1,\ldots, S_d)$ be a random vector with joint survival function $\bar{G}_d$ as in \eqref{eq:Ssurvivalfunction} for some stdf $\ell$. Then $1/\bm{S}$ belongs to the maximum domain of attraction of a multivariate extreme-value distribution with unit Fr\'echet margins and stdf $\ell$.
\end{proposition}
\begin{proof}[Proof of Proposition~\ref{prop:stdfS}]
	For the margins, recall that for each $i \in \{1,\ldots, d\}$, $S_i \sim \mathrm{Beta}(1,d-1)$. The survival function of $1/S_i$ is thus given by
		$\bar{F}_{1/S_i}(s)=1-(1-1/s)^{d-1}$; it is easily seen that {$\bar{F}_{1/S_i}\in\mathcal{R}_{-1}$}. Now set $c_n=\{1-(1-1/n)^{1/(d-1)} \}^{-1}$. From Equation 3.13 in \cite{Embrechts/Kluppelberg/Mikosch:1997}, for all $s_i \in \mathbb{R}$, it then holds that $\Pr(1/S_i \le c_n s_i) \to \Phi_1(s_i)$ as $n \to \infty$. Thus $1/
	\bm{S}$ is in the domain of attraction of a multivariate extreme-value distribution with unit Fr\'echet margins and stdf $\ell$ if and only if for all $\bm{s} \in \mathbb{R}_+^d$,
	\begin{multline*}
	\lim\limits_{n\to\infty}n\left\{1-\Pr(1/S_1\le c_ns_1,\ldots,1/S_d\le c_ns_d)\right\} \\= \lim\limits_{n\to\infty}n \left[1-\bar G_d\bigl\{1/(c_n s_1), \ldots, 1/(c_n s_d)\bigr\}\right] = \ell(1/s_1,\ldots,1/s_d).
	\end{multline*} 
	To show this, fix an arbitrary $\bm{s} \in \mathbb{R}_+^d$ and observe that because $c_n \to \infty$ as $n \to \infty$, 
	$$
	\bar G_d\bigl\{1/(c_n s_1), \ldots, 1/(c_n s_d)\bigr\} = \bigl\{1- (1/c_n) \ell(1/s_1,\ldots, 1/s_d)\bigr\}^{d-1}
	$$
	for all $n$ sufficiently large. Now note that as $n \to \infty$, $n/c_n^k$ converges to $0$ for all $k \in \{2,\ldots, d-1\}$ and to $1/(d-1)$ for $k=1$. Consequently,
	\begin{multline*}
	\lim\limits_{n\to\infty}n \left[ 1-  \bigl\{1- (1/c_n) \ell(1/s_1,\ldots, 1/s_d)\bigr\}^{d-1} \right] \\ = \lim\limits_{n\to\infty} \sum_{k=1}^{d-1} {d-1 \choose k} (-1)^{k+1} \frac{n}{c_n^k} \ell^k(1/s_1,\ldots, 1/s_d) =  \ell(1/s_1,\ldots, 1/s_d)
	\end{multline*}
	as claimed.
\end{proof}
The following lemma determines the normalizing sequences needed for the proof of Theorem~\ref{thm:MaxDofA}.
\begin{lemma}\label{lemma:marginal}
	Let $C_{\mathcal{G},\bm{\psi},\bm{\ell},Q}$ be a clustered Archimax copula such that Assumptions~\ref{assumption:MDA} and \ref{assumption:stdf} are satisfied. Then the following hold:
	\begin{itemize}
		\item[(i)] For each $k\in\mathcal{D}_1$ and $i\in\{1,\ldots,d_k\}$, $1/(R_kS_{ki})\in\mathcal{M}(\Phi_{\rho_k})$. Recall that for $k\in\mathcal{D}_1$, $b_k= \E\{1/Z_k^{\rho_k} \}$ where $Z_k\sim \mathrm{Beta}(1,d_k-1)$. Moreover, there exists a sequence of positive constants $\{a_{nk}\}$ such that for all $x>0$, $n\Pr(1/R_k>a_{nk}x)\rightarrow x^{-\rho_k}$ as $n\to\infty$ and $n\Pr(1/(R_kS_{ki})>a_{nk}b_k^{1/\rho_k}x)\rightarrow x^{-\rho_k}$ as $n\to\infty$.
		\item[(ii)]  For each $k\in\mathcal{D}_2$ and $i\in\{1,\ldots,d_k\}$, $1/(R_kS_{ki})\in\mathcal{M}(\Phi_{1})$. Moreover, there exists a sequence of positive constants $\{a_{nk}\}$ such that for all $x>0$, $n\Pr(1/S_{ki}>a_{nk}x)\rightarrow x^{-1}$ as $n\to\infty$ and $n\Pr(1/(R_kS_{ki})>a_{nk}{b_k}x)\rightarrow x^{-1}$ as $n\to\infty$, where $b_k=\E(1/R_k)$.
	\end{itemize}
\end{lemma}
\begin{proof}[Proof of Lemma~\ref{lemma:marginal}]
	\noindent(i)	Let $k\in\mathcal{D}_1$ and $i\in\{1,\ldots,d_k\}$. We then have $(1/R_k)\in\mathcal{M}(\Phi_{\rho_k})$ by assumption and $1/S_{ki}\in\mathcal{M}(\Phi_{1})$ owing to the fact that $S_{ki} \sim \mathrm{Beta}(1,d-1)$. By Proposition 3.1.1 in \cite{Embrechts/Kluppelberg/Mikosch:1997}, there exists a sequence of positive constants $\{a_{nk}\}$ such that for all $x>0$, $n\Pr(1/R_k>a_{nk}x)\rightarrow x^{-\rho_k}$ as $n\to\infty$. Because $\rho_k < 1$, $\E(1/S_{ki})^{\rho_k + \varepsilon} < \infty$ for some $\varepsilon$ sufficiently small. Using the lemma of \cite{Breiman:1965} and the fact that $b_k=\E \{1/S_{ki}^{\rho_k} \}$, we then have, for all $x > 0$ and with $\zeta_{kn} = a_{nk}b_k^{1/\rho_k}$,
	\begin{align*} %\label{eq:marginsD1}
	\lim\limits_{n\rightarrow\infty}n\Pr\Bigl(\frac{1}{R_kS_{ki}} > \zeta_{kn}x\Bigr) &= \lim\limits_{n\rightarrow\infty}n\Pr\Bigl(\frac{1}{R_k} > \zeta_{kn}x\Bigr) \frac{\Pr\Bigl(\frac{1}{R_kS_{ki}} > \zeta_{kn}x \Bigr)}{\Pr\Bigl(\frac{1}{R_k} > \zeta_{kn}x \Bigr)}\\
	&=(xb_k^{1/\rho_k})^{-\rho_k}b_k = x^{-\rho_k}\;.
	\end{align*}
Indeed, $n\Pr(1/R_k > a_{nk}b_k^{1/\rho_k}x) \to (xb_k^{1/\rho_k})^{-\rho_k}$ as $n\to\infty$ by the choice of normalizing constants $\{a_{nk}\}$. The convergence of the fraction in the above display is due to Breiman's lemma. The Fisher-Tippett-Gnedenko Theorem \citep{Fisher/Tippett:1928,Gnedenko:1943} implies that since $1/R_k\in\mathcal{M}({\Phi}_{\rho_k})$ and $\rho_k\in(0,1)$, $\bar{F}_{1/R_k}\in\mathcal{R}_{-\rho_k}$. We also have that $1/S_{ki}$ and ${1/R_k}$ are independent, positive, and $\E(1/S_{ki}^\gamma] < \infty $ for $\gamma\in(\rho_k,1)$. By Breiman's lemma, $1/(R_k S_{ki}) \in \mathcal{M}(\Phi_{\rho_k})$ and 
		$$
		\frac{\Pr\Bigl(\frac{1}{R_kS_{ki}} > a_{nk}b_k^{1/\rho_k}x \Bigr)}{\Pr\Bigl(\frac{1}{R_k} > a_{nk}b_k^{1/\rho_k}x \Bigr)} \to \E(1/S_{ki}^{\rho_k}) = b_k
		$$
		as $n\to\infty$.	

	\medskip	
	\noindent(ii) 	Let $k\in\mathcal{D}_2$ and $i\in\{1,\ldots,d_k\}$. The proof of the result relies again on  Breiman's lemma; see also Proposition 2(b) of \cite{Belzile/Neslehova:2017}. 
	Since $1/S_{ki} \in \mathcal{M}(\Phi_{1})$, Proposition 3.1.1 in \cite{Embrechts/Kluppelberg/Mikosch:1997} implies that there exist sequences of positive constants $\{a_{nk}\}$ such that for all $x>0$, $n\Pr(1/S_{ki} > a_{nk}x) \rightarrow x^{-1}$ as $n\to\infty$, and this for all $i=1,\ldots,d_k$. Recall that $b_k = \E(1/R_k)$. Similarly to the proof of part (i), Breiman's lemma then implies that for all $x > 0$ and with $\zeta_{kn} = a_{nk}b_k$,
	\begin{align*} %\label{eq:marginsD2}
	\lim_{n\rightarrow\infty} n\Pr\Bigl(\frac{1}{R_kS_{ki}} > \zeta_{kn}x\Bigr)
	&= \lim_{n\rightarrow\infty} n\Pr\Bigl(\frac{1}{S_{ki}} > \zeta_{kn}x \Bigr)\frac{\Pr\Bigl(\frac{1}{R_kS_{ki}} > \zeta_{kn}x \Bigr)}{\Pr\Bigl(\frac{1}{S_{ki}} > \zeta_{kn}x \Bigr)}\\
	&= (xb_k)^{-1}b_k = x^{-1}\;.
	\end{align*}
The convergence of the first part of the above is due to the choice of the normalizing constants $\{a_{nk}\}$. For the convergence of the second term, note that $\bar{F}_{1/S_{ki}} \in \mathcal{R}_{-1}$ and by assumption, $\E\{1/R_k^{1+\epsilon_k} \}$ for some $\epsilon_k>0$. Finally, since $1/S_{ki}$ and $1/R_k$ are independent and positive, Breiman's lemma implies that $1/(R_k S_{ki}) \in \mathcal{M}(\Phi_1)$ and that
		$$
		\frac{\Pr\Bigl(\frac{1}{R_kS_{ki}} > a_{nk}b_kx\Bigr)}{\Pr\Bigl(\frac{1}{S_{ki}} > a_{nk}b_kx\Bigr)} \to \E (1/R_k) = b_k
		$$
		as $n\to\infty$. This completes the proof.
\end{proof}
The lemma below establishes asymptotic independence between clusters in $\mathcal{D}_1$ and clusters in $\mathcal{D}_2$.
\begin{lemma}\label{lemma:D1D2}
	Suppose that $k\in\mathcal{D}_1$, $l\in\mathcal{D}_2$, $i\in \{1,\ldots,d_k\}$ and $j\in\{1,\ldots,d_l \}$. Let $\{a_{nk}\}$ and $\{a_{nj}\}$ be normalizing sequences as in Lemma~\ref{lemma:marginal}. As in Lemma~\ref{lemma:marginal} (ii), let $b_l=\E(1/R_l) $. Then for all $x,y>0$,
	$$
	\lim_{n\to\infty}n\Pr\{1/(R_kS_{ki}) > a_{nk}b_k^{1/\rho_k}x,\ 1/(R_lS_{lj}) > a_{nl}b_l y\} = 0 \;.
	$$
\end{lemma}
\begin{proof}
	Fix $x, y > 0$ and recall that $\rho_k\in(0,1)$. The probability of interest can be written as follows, for $\zeta_{nk} = a_{nk}b_k^{1/\rho_k}$ and $\eta_{nl} = a_{nl}b_l$,
	\begin{align*}
	& n\Pr\left\{1/(R_kS_{ki}) > \zeta_{nk}x,\ 1/(R_lS_{lj}) > \eta_{nl}y \right\}\\
	& = \int_{\mathbb{R}_+^2} n\Pr\left(1/S_{ki} > \zeta_{nk}x r_k,\ 1/S_{lj} > \eta_{nl}y r_l\right) dF_{R_k,R_l}(r_k,r_l)\\
	& =\int_{\mathbb{R}_+^2} n\Pr\{1/S_{ki} > \zeta_{nk}x r_k\}\Pr\{ 1/S_{lj}> \eta_{nl}y r_l\}dF_{R_k,R_l}(r_k,r_l)\;,
	\end{align*}
	where the first equality is due to the independence between $(R_k,R_l)$ and $(S_{ki},S_{lj})$ and the last equality is due to the independence of $S_{ki}$ and $S_{lj}$. Next, consider the integrand as a sequence of functions $\{f_n\}$ defined on $\mathbb{R}_+^2$. Observe that for each $r_k, r_l > 0$,
	$$
	f_n(r_k,r_l)\le g_n(r_k,r_l)\;,
	$$
	where $\{g_n\}$ is itself a sequence of functions on $\mathbb{R}_+^2$ defined by
	$$
	g_n(r_k,r_l)=g_n(r_l)=n\Pr\{ 1/S_{lj} > a_{nl}b_lyr_l\}\;.
	$$
	From the choice of $\{a_{nl}\}_{n\in\mathbb{N}}$, for all $r_k,r_l>0$,
	$\lim_{n\to\infty}g_n(r_k,r_l)=g(r_k,r_l)$, where $g(r_k,r_l)=1/(b_lyr_l)$. Moreover, 
	$$\int_{\mathbb{R}_+^2}g(r_k,r_l)dF_{R_k,R_l}(r_k,r_l)=\int_{\mathbb{R}_+^2}\frac{1}{b_lyr_l}dF_{R_k,R_l}(r_k,r_l)=\frac{1}{y}\;,$$
	and
	$$\int_{\mathbb{R}_+^2}g_n(r_k,r_l)dF_{R_k,R_l}(r_k,r_l)=n\Pr\{ 1/(R_lS_{lj})>a_{nl}b_l y\}\to \frac{1}{y}
	$$
	as $n\to\infty$. We therefore have a sequence of nonnegative functions $\{g_n\}$ bounding $\{f_n\}$ from above such that 
	$$
	\lim\limits_{n\to\infty}\int_{\mathbb{R}_+^2}g_n(r_k,r_l)dF_{R_k,R_l}(r_k,r_l)=\int_{\mathbb{R}_+^2}\lim\limits_{n\to\infty}g_n(r_k,r_l)dF_{R_k,R_l}(r_k,r_l)\;.
	$$ 
	Finally, note that
	\begin{equation*}
	f_n(r_k,r_l)=n\Pr\{1/S_{ki} > a_{nk}b_k^{1/\rho_k}x r_k \} \Pr\{ 1/S_{lj}> a_{nl}b_l y r_l\} \to 0
	\end{equation*}
	as $n\to\infty$ since 
	\begin{align*}
	\Pr\{1/S_{ki} > a_{nk}b_k^{1/\rho_k}xr_k \} \to 0 \quad \text{and} \quad
	n\Pr\{ 1/S_{lj} > a_{nl}b_lyr_l \} \to 1/(b_lyr_l)
	\end{align*}
	as $n\to\infty$. The desired result then follows by the generalized Lebesgue dominated convergence theorem (see Theorem 1.21 in \cite{Kallenberg:2002}, for example).
\end{proof}
We are now ready to prove Theorem~\ref{thm:MaxDofA}.

\subsection{Proof of Theorem \ref{thm:MaxDofA}} \label{sec:proofmda}
A random vector $(Y_1,\ldots,Y_d)$ is in the maximum domain of attraction of the extreme-value distribution $H$ with Fr\'echet margins if and only if there exist sequences of positive constants $(a_{ni})\in(0,\infty)$, $i\in\{1,\ldots,d\}$, so that, for all $(y_1,\ldots,y_d)\in\mathbb{R}_+^d$,
	\begin{equation*}
	\lim\limits_{n\rightarrow\infty}n\left\{1-\Pr\left(Y_1\le a_{n1}y_1,\ldots,Y_d\le a_{nd}y_d\right)\right\}=-\ln H(y_1,\ldots,y_d)\;.
	\end{equation*}
This is a multivariate extension of Proposition 3.1.1 in \cite{Embrechts/Kluppelberg/Mikosch:1997}, as used in \cite{Belzile/Neslehova:2017}. For each $k\in\{1\,\ldots,K \}$, set the sequences $\{a_{nk}\}$ as done in Lemma~\ref{lemma:marginal}.
Then the fact that the marginals of $H$ are Fr\'echet follows from the said Lemma. 
With the normalizing constants now set, the limit of interest is, for any fixed $\bm{x} \in \mathbb{R}_+^d$ and, for convenience, let $\mathcal{J}=\{(k,i):k=1,\ldots,K, i=1\ldots,d_k \}$,
\begin{equation} \label{eq:fulllimit}
\lim_{n\to\infty} n\left[ 1 - \Pr\left( \bigcap_{(k,i) \in \mathcal{J}} \left\{ \frac{1}{R_kS_{ki}} \le a_{nk} b_k^{1/\rho_k}x_{ki}\right\} \right) \right]\;,
\end{equation}
where for $k\in\mathcal{D}_2$, $b_k=\E(1/R_k) $ as in Lemma~\ref{lemma:marginal} (ii) and for ease of notation, $\rho_k=1$.
Letting $\mathcal{P}(\mathcal{J})$ denote the power set of $\mathcal{J}$, \eqref{eq:fulllimit} can be rewritten {as}
\begin{equation}\label{eq:bigsum}
\lim_{n\to\infty} n\sum_{p\in\mathcal{P}(\mathcal{J})} (-1)^{|p|+1}\Pr\left(\bigcap_{{(k,i)}\in p} \left\{\frac{1}{R_kS_{ki}} > a_{nk}b_k^{1/\rho_k}x_{ki}\right\} \right)\;.
\end{equation}
We now consider the subset $\mathcal{P}({\mathcal{J}})|_{\mathcal{D}_1,\mathcal{D}_2}$ of $\mathcal{P}({\mathcal{J}})$ such that $p \in \mathcal{P}(\mathcal{J})|_{\mathcal{D}_1,\mathcal{D}_2}$ if and only if there exists $(k,i),(l,j)\in p$ so that $k\in\mathcal{D}_1$ and $l\in\mathcal{D}_2$.
As it turns out, the summands in \eqref{eq:bigsum} for which $p\in \mathcal{P}(\mathcal{J})|_{\mathcal{D}_1,\mathcal{D}_2}$ are asymptotically negligible.
To see this, fix an arbitrary $p\in \mathcal{P}(\mathcal{J})|_{\mathcal{D}_1,\mathcal{D}_2}$ and pick $(k,i),(l,j) \in p$ so that $k\in\mathcal{D}_1$ and $l\in\mathcal{D}_2$. Then for all $\{x_{ca}:{(c,a)}\in p\}\in\mathbb{R}_+^{|p|}$,
\begin{align*}
n\Pr&\left(\bigcap_{(c,a)\in p} \left\{\frac{1}{R_cS_{ca}} > a_{nc} b_c^{1/\rho_c}x_{ca}\right\} \right)\\
&\le n\Pr\{1/(R_kS_{ki})>a_{nk}b_k^{1/\rho_k}x_{ki}, 1/(R_lS_{lj})>a_{nl}b_l^{1/\rho_l}x_{lj}\}\to 0
\end{align*}
as $n\to\infty$ by Lemma \ref{lemma:D1D2}.

Now, let $\mathcal{P}(\mathcal{J})|_{\mathcal{D}_1}$ be the subset of $\mathcal{P}(\mathcal{J})$ such that $p\in \mathcal{P}(\mathcal{J})|_{\mathcal{D}_1}$ if and only if $(c,a)\in p$ implies that $c\in\mathcal{D}_1$.
In other words, $\mathcal{P}(\mathcal{J})|_{\mathcal{D}_1}$ contains only sets of indices $(c,a)$ with $c\in\mathcal{D}_1$. Let $N_1=\sum_{k\in\mathcal{D}_1}d_k$, define $\bm{S}^{(1)} = (\bm{S}_k : k \in \mathcal{D}_1)$, and rewrite the summands in \eqref{eq:bigsum} with $p\in \mathcal{P}(\mathcal{J})|_{\mathcal{D}_1}$ as follows:
\begin{align*}
&n\sum_{p\in\mathcal{P}(\mathcal{J})|_{\mathcal{D}_1}} (-1)^{|p|+1} \Pr\Biggl(\bigcap_{(c,a)\in p} \Biggl\{ \frac{1}{R_cS_{ca}} > a_{nc}b_c^{1/\rho_c}x_{ca}\Biggr\} \Biggr)\\
&= n\Biggl[ 1 - \Pr\Biggl(\bigcap_{k\in\mathcal{D}_1} \bigcap_{1 \leq i \leq d_k} \Biggl\{ \frac{1}{R_kS_{ki}} \le a_{nk}b_k^{1/\rho_k}x_{ki} \Biggr\} \Biggr) \Biggr]\\
&=\int_{[0,1]^{N_1}}n\Biggl[ 1 - \Pr\Biggl(\bigcap_{k\in\mathcal{D}_1} \bigcap_{1 \leq i \leq d_k} \Biggl\{ \frac{1}{R_k} \le a_{nk}b_k^{1/\rho_k}x_{ki}s_{ki} \Biggr\} \Biggr) \Biggr] dF_{\bm{S}^{(1)}}(\bm{s}^{(1)})\\
&=\int_{[0,1]^{N_1}}n\Biggl[ 1 - \Pr\Biggl(\bigcap_{k\in\mathcal{D}_1} \Biggl\{ \frac{1}{R_k} \le a_{nk}b_k^{1/\rho_k} \min_{1 \leq i \leq d_k} (x_{ki}s_{ki}) \Biggr\} \Biggr) \Biggr] dF_{\bm{S}^{(1)}}(\bm{s}^{(1)})\;,
\end{align*}
where $\bm{s}^{(1)} \in [0,1]^{N_1}$ is indexed so that $\bm{s}^{(1)}=(s_{ki}: k \in \mathcal{D}_1,\ 1 \leq i \leq d_k)$.
Now consider the integrand as a sequence of functions $\{f_n\}$ defined on $[0,1]^{N_1}$ and observe that for each $n \in \mathbb{N}$, $0 \le f_n \le g_n$, where $g_n$ is given, for each $\bm{s}^{(1)}\in[0,1]^{N_1}$ as above, by
$$
g_n(\bm{s}^{(1)})= n \sum_{k\in\mathcal{D}_1} \sum_{i=1}^{d_k} \Pr\left(1/R_k> a_{nk}b_k^{1/\rho_k}x_{ki} s_{ki} \right) 
$$
Clearly,	 $g_n(\bm{s}^{(1)})\to g(\bm{s}^{(1)})$ as $n\to\infty$, where 
$$
g(\bm{s}^{(1)}) = \sum_{k\in\mathcal{D}_1}\sum_{i=1}^{d_k} \frac{1}{b_k\{x_{ki} s_{ki}\}^{\rho_k}}
$$
with
\begin{align*}
\int_{[0,1]^{N_1}} g(\bm{s}^{(1)})dF_{\bm{S}^{(1)}}(\bm{s}^{(1)})
=\sum_{k\in\mathcal{D}_1}\sum_{i=1}^{d_k} x_{ki}^{-\rho_k}\;.
\end{align*}
Moreover,
\begin{align*}
\int_{[0,1]^{N_1}}g_n(\bm{s}^{(1)}) dF_{\bm{S}^{(1)}}(\bm{s}^{(1)})
&= n\left\{\sum_{k\in\mathcal{D}_1}\sum_{i=1}^{d_k}\Pr\left(1/(R_k S_{ki} ) > a_{nk}b_k^{1/\rho_k}x_{ki}\right) \right\}\\
&\to \sum_{k\in\mathcal{D}_1}\sum_{i=1}^{d_k} x_{ki}^{-\rho_k}
\end{align*}
as $n\to\infty$. Therefore, we have a sequence of majorants $\{g_n\}$ for which the identity $\lim_{n\to\infty}\int g_n =\int \lim_{n\to\infty} g_n$ is satisfied.
Now recall that the vector of reciprocal distortions $1/\bm{R}$ has a limiting stdf $\ell_{1/\bm{R}}$ defined in terms of the positive, unit-mean variables $W_1,\ldots, W_K$ in Assumption~\ref{assumption:stdf}.
Therefore, $f_n \to f$ pointwise, where for all $\bm{s}^{(1)} \in[0,1]^{N_1}$,
\begin{align*}
f(\bm{s}^{(1)})&= \E\biggl\{ \max_{k\in\mathcal{D}_1}\frac{W_k}{b_k \min\limits_{1 \leq i \leq d_k}(x_{ki} s_{ki})^{\rho_k}} \biggr\} \;.
\end{align*}
Now, integrating over the $\bm{s}^{(1)}$ yields
\begin{equation*}
\int_{[0,1]^{N_1}} f(\bm{s}^{(1)}) dF_{\bm{S}^{(1)}}(\bm{s}^{(1)})=\E\Biggl\{ \max_{\substack{k\in\mathcal{D}_1 \\ 1 \leq i \leq d_k}} \frac{W_k}{b_k (S_{ki} x_{ki})^{\rho_k}} \Biggr\} \;.
\end{equation*}
Using the generalized Lebesgue dominated convergence theorem, we can {thus} conclude that for all $\bm{x}^{(1)} = (\bm{x}_k:k\in\mathcal{D}_1 )\in\mathbb{R}_+^{N_1}$,
\begin{align*}
\lim_{n\to\infty} n\sum_{p\in\mathcal{P}(\mathcal{J})|_{\mathcal{D}_1}} (-1)^{|p|+1} \Pr\Bigl(\bigcap\limits_{(c,a) \in p} \{1/(R_cS_{ca}) > a_{nb}b_c^{\rho_c}x_{ca}\} \Bigr)\\
= \E\Biggl\{ \max_{\substack{k\in\mathcal{D}_1\\ 1 \leq i \leq d_k}} \frac{W_k}{b_k(S_{ki} x_{ki})^{\rho_k}} \Biggr\}\;.
\end{align*}

Analogously to $\mathcal{P}(\mathcal{I})|_{\mathcal{D}_1}$, let $\mathcal{P}(\mathcal{I})|_{\mathcal{D}_2}$ contain only sets of indices $(c,a)$ with $c\in\mathcal{D}_2$.
Let $K_2=|\mathcal{D}_2|$ and $N_2=\sum_{k\in\mathcal{D}_2}d_k$, define $\bm{R}^{(2)} = (R_k : k \in \mathcal{D}_2)$, and recall that $\rho_k=1$ for all $k\in\mathcal{D}_2$.
Next, for $\bm{x}^{(2)} = (\bm{x}_k : k\in\mathcal{D}_2) \in \mathbb{R}_+^{N_2}$, rewrite the summands of \eqref{eq:bigsum} with $p \in \mathcal{P}(\mathcal{I})|_{\mathcal{D}_2}$ as follows:
\begin{align*}
& n \sum_{p\in\mathcal{P}(\mathcal{I})|_{\mathcal{D}_2}} (-1)^{|p|+1}\Pr\Bigl(\bigcap_{(c,a)\in p} \{1/(R_cS_{ca}) > a_{nc}b_c^{1/\rho_c}x_{ca}\} \Bigr)\\
&=n\biggl[1-\Pr\biggl(\bigcap\limits_{k\in\mathcal{D}_2} \bigcap\limits_{1 \leq i \leq d_k} \{1/(R_kS_{ki}) \le a_{nk}b_kx_{ki} \biggr)\biggr]\\
&=\int_{\mathbb{R}_+^{K_2}}n\biggl[1-\Pr\biggl(\bigcap_{k\in\mathcal{D}_2}\bigcap_{1 \leq i \leq d_k} \{1/S_{ki} {\le} a_{nk}b_kx_{ki}r_k\} \biggr)\biggr] dF_{\bm{R}^{(2)}}(\bm{r})\;,
\end{align*}
where $\bm{r}^{(2)} = (r_k: k \in \mathcal{D}_2)$ is indexed as $\bm{R}^{(2)}$.
Now consider the integrand as a sequence of functions $\{f_n\}$ defined on $\mathbb{R}_+^{K_2}$ and observe that for each $n \in \mathbb{N}$, $0 \le f_n \le g_n$, where $g_n$ is given, for all $\bm{r}^{(2)} \in\mathbb{R}_+^{K_2}$, by
$$
g_n(\bm{r}^{(2)})=n\Bigl\{\sum_{k\in\mathcal{D}_2}\sum_{i=1}^{d_k}\Pr\bigl(1/S_{ki} > a_{nk}b_k x_{ki} r_k \bigr) \Bigr\}
$$
Clearly, for all $\bm{r}^{(2)} \in \mathbb{R}_+^{K_2}$ and as $n \to \infty$,
$$
g_n(\bm{r}^{(2)}) \to g(\bm{r}^{(2)}) = \sum_{k\in\mathcal{D}_2}\sum_{i=1}^{d_k} \frac{1}{b_kx_{ki} r_k} \;.
$$
Furthermore,
\begin{align*}
\int_{\mathbb{R}_+^{K_2}}g(\bm{r}^{(2)})dF_{\bm{R}^{(2)}}(\bm{r}^{(2)}) = \sum_{k\in\mathcal{D}_2}\sum_{i=1}^{d_k} \frac{1}{x_{ki}}
\end{align*}
and, as $n \to \infty$,
\begin{align*}
\int_{\mathbb{R}_+^{K_2}} g_n(\bm{r}^{(2)})
dF_{\bm{R}^{(2)}}(\bm{r}^{(2)}) &= n\Bigl\{\sum_{k\in\mathcal{D}_2}\sum_{i=1}^{d_k} \Pr\bigl(1/(R_kS_{ki}) > a_{nk}b_kx_{ki}\bigr) \Bigr\}\\
	&\to \sum_{k\in\mathcal{D}_1}\sum_{i=1}^{d_k} \frac{1}{x_{ki}}\;.
\end{align*}
Analogously to the treatment of $\mathcal{P}(\mathcal{I})|_{\mathcal{D}_1}$, we have a sequence of majorants $\{g_n\}$ such that $\lim_{n\to\infty}\int g_n =\int \lim_{n\to\infty} g_n$. It remains to determine the limit of the sequence of functions $\{f_n\}$ defined for all $\bm{r}^{(2)} \in \mathbb{R}_+^{K_2}$ by
$$
n\biggl[1-\Pr\biggl(\bigcap_{k\in\mathcal{D}_2}\bigcap_{1 \leq i \leq d_k} \{1/S_{ki} > a_{nk} b_k x_{ki} r_k\} \biggr)\biggr] \;.
$$
By assumption, $\bm{S}_k$ and $\bm{S}_l$ are independent if $k\ne l$ and are therefore asymptotically independent as well.
Using Proposition \ref{prop:stdfS} and the fact that $1/S_{ki}\in\mathcal{M}(\Phi_1)$ for all $k\in\{1,\ldots,K\}$ and $i\in\{1,\ldots,d_k\}$, one has that $f_n \to f$ pointwise, where for all $\bm{r}^{(2)} \in \mathbb{R}_+^{K_2}$,
\begin{align*}
f(\bm{r}^{(2)})&= \sum_{k\in\mathcal{D}_2} \ell_k\left\{(b_kx_{k1}r_k)^{-1},\ldots,(b_kx_{kd_k}r_k)^{-1}\right\}\\
&=\sum_{k\in\mathcal{D}_2}(b_k r_k)^{-1} \ell_k\left(x_{k1}^{-1},\ldots,x_{kd_k}^{-1}\right)\;.
\end{align*}
Integrating the limit $f$ yields
\begin{equation*}
\int_{\mathbb{R}_+^{K_2}}f(\bm{r}^{(2)})
dF_{\bm{R}^{(2)}}(\bm{r}^{(2)}) = \sum_{k\in\mathcal{D}_2} \ell_k\left(x_{k1}^{-1},\ldots,x_{kd_k}^{-1}\right)\;.
\end{equation*}
Thus for all $\bm{x} \in \mathbb{R}_+^d$, the limit  \eqref{eq:fulllimit} is equal to
$$
\E\Biggl\{ \max_{\substack{k\in\mathcal{D}_1 \\ 1 \leq i \leq d_k}} \frac{W_k}{b_k (S_{ki} x_{ki})^{\rho_k}} \Biggr\} + \sum\limits_{k\in\mathcal{D}_2}\ell_k\left(x_{k1}^{-1},\ldots,x_{kd_k}^{-1}\right)\;.
$$
Recalling that $1/(R_kS_{ki})\in\mathcal{M}(\Phi_{\rho_k})$, one obtains \eqref{eq:stdf} by plugging in the appropriate Fr\'echet margins.
\qed

\subsection{Proofs of Corollaries \ref{cor:AI} and \ref{cor:asymptoticordering}}\label{sec:sec:proofcor}
\begin{proof}[Proof of Corollary \ref{cor:AI}]
	Let $K_1=|\mathcal{D}_1|$ and recall that $\E[\max_{k\in\mathcal{D}_1} y_kW_k]$ is the limiting stdf of $\{1/R_k:k\in\mathcal{D}_1\}$, defined for all $(y_1,\ldots,y_{K_1})\in\mathbb{R}_+^{K_1}$. Letting $\{W_k:k\in\mathcal{D}_1\}$ be a (uniformly) random permutation of $(K_1,0,\ldots,0)$ yields the independence stdf $\E[\max_{k\in\mathcal{D}_1} y_kW_k]=y_1+\ldots+y_{K_1}$. Due to the fact that, for all $k,l \in \{1,\dots,d\}$ and $i \in \{1,\dots,d_l\}$, $W_k$ is independent of $S_{li}$, plugging this into \eqref{eq:stdf} yields, for all $\bm{x}\in\mathbb{R}_+^d$,
	\begin{align*}
	\ell_{\mathcal{G},\bm{\psi},\bm{\ell},Q}(\bm{x})& =\sum_{k\in\mathcal{D}_1}\E\Bigl( \max_{i=1,\ldots,d_k} \frac{x_{ki}}{b_kS_{ki}^{\rho_k}}\Bigr) + \sum_{k\in\mathcal{D}_2} \ell_k(\bm{x}_k)\;.
	\end{align*}
	By Proposition~6.1 in \cite{Charpentier/Fougeres/Genest/Neslehova:2014}, as mentioned in the Section~\ref{sec:intro}, for each $k\in\mathcal{D}_1$,
	$$
	\E\Bigl( \max_{1 \leq i \leq d_k} \frac{x_{ki}}{b_kS_{ki}^{\rho_k}}\Bigr) = \ell_k^{\rho_k}(\bm{x}_k^{1/\rho_k})\;.
	$$
	{This completes the proof.}\end{proof}

\begin{proof}[Proof of Corollary \ref{cor:asymptoticordering}]
Observe first that by Assumption \ref{assumption:stdf} and the fact that the variables $W_k$ have unit mean,
	\begin{align}\label{eq:ellRuseful}
	\ell_{1/\bm{R}}(\bm{x}') &= \E\bigl\{\max_{1 \leq k \leq K}(x_k' W_k)\bigr\} \nonumber\\
	 &\le \E\bigl\{\max_{k\in \mathcal{D}_1}(x_k' W_k) + \sum_{k\in \mathcal{D}_2} x_k' W_k\bigr\} \nonumber\\
	 &= \E\bigl\{\max_{k\in \mathcal{D}_1}(x_k' W_k)\bigr\} + \sum_{k\in \mathcal{D}_2} x_k'.
	\end{align}
	Next, note that $\ell_{\mathcal{G},\bm{\psi},\bm{\ell},Q}(\bm{x}) = \mathcal{A}(\bm{x}') + \mathcal{B}(\bm{x}')$, where
	$$
	\mathcal{A}(\bm{x}') = \sum_{k \in \mathcal{D}_2} \ell_k \bigl(\bm{x}_k) = \sum_{k \in \mathcal{D}_2} x_{i_k} = \sum_{k \in \mathcal{D}_2} x_{k}'
	$$
	and
	$$
	\mathcal{B}(\bm{x}) = \E\Biggl( \max_{\substack{k  \in \mathcal{D}_1 \\ 1 \leq j \leq d_k}} \frac{\bm{x}_{kj} W_k}{b_k S_{kj}^{\rho_k}} \Biggr) = \E\biggl( \max_{k  \in \mathcal{D}_1}  \frac{x_k' W_k}{b_k S_{kj}^{\rho_k}} \biggr).
	$$
	Because for each $k \in \mathcal{D}_1$, $b_k = \E(1/S_{kj}^{\rho_k})$, we have that for any $\bm{w} \in \mathbb{R}_+^K$ and $k \in \mathcal{D}_1$,
	$$
	\E\biggl( \max_{k  \in \mathcal{D}_1}  \frac{x_k' w_k}{b_k S_{kj}^{\rho_k}} \biggr) \ge \E\biggl(\frac{x_k' w_k}{b_k S_{kj}^{\rho_k}} \biggr) = x_k' w_k \;,
	$$
	so that 
	$$
	\E\biggl( \max_{k  \in \mathcal{D}_1}  \frac{x_k' w_k}{b_k S_{kj}^{\rho_k}} \biggr) \ge \max_{k \in \mathcal{D}_1} (x_k' w_k)\;.
	$$
	This implies that 
	$$
	\mathcal{B}(\bm{x}') \ge  \E\Bigl\{\max_{k\in \mathcal{D}_1}(x_k' W_k)\Bigr\}\;,
	$$
	which together with \eqref{eq:ellRuseful} yields the desired result.
\end{proof}

\section{Jackknife estimation of the covariance matrix of $\bm{T}$} \label{app:jackknife}

Here, we provide a consistent estimator of the asymptotic variance $\bm{\Sigma}$ of $\sqrt{n} \bm{T}$, where $\bm{T}$ is as in \eqref{eq:T}.
Recall the definition of $\mathcal{I}$ in \eqref{eq:I}.
For each $\nu \in \{1,\dots,n\}$ and $(i,j,k) \in \mathcal{I}$, let $\hat{\theta}_{ij\nu}$ be a version of $\hat{\theta}_{ij}$ based on all but the $\nu$th observation; $\bar{\theta}_{k\nu}$ be the average of all $\hat{\theta}_{ij\nu}$ such that, given $k$, $(i,j,k) \in \mathcal{I}$; $T_{ijk\nu} = \hat{\theta}_{ij\nu} - \bar{\theta}_{k\nu}$; and $\bm{T}_{\nu} = (T_{\iotaa\nu})_{\iotaa \in \mathcal{I}}$.
In particular, $\bm{T}_{\nu}$ is a version of $\bm{T}$ based on all but the $\nu$th observation.
It then follows from Theorem~9 of \cite{Arvesen:1969} that, for $\bm{T}_{\nu}^* = n \bm{T} - (n-1) \bm{T}_{\nu}$ and $\bm{T}_{\bullet}^* = (1/n) \sum_{\nu=1}^n \bm{T}_{\nu}^*$,
\begin{equation} \label{eq:T-covariance}
	\hat{\bm{\Sigma}} = \frac{1}{n-1} \sum_{\nu=1}^n (\bm{T}_{\nu}^* - \bm{T}_{\bullet}^*)(\bm{T}_{\nu}^* - \bm{T}_{\bullet}^*)^\top
\end{equation}
is a consistent estimator of $\bm{\Sigma}$.

\section{Conjectured extension of Theorem~\ref{thm:MaxDofA}}\label{sec:conjecture}
As it is stated, Theorem~\ref{thm:MaxDofA} does not account for the boundary case when $1/R_k \in\mathcal{M}(\Phi_{1})$, which can occur. It would thus be desirable to replace Assumption~\ref{assumption:MDA} of Theorem~\ref{thm:MaxDofA} by the following requirement.
	\begin{assumption}\label{assumption:MDAExt}
		For a clustered Archimax copula as in Definition~\ref{def:model}, assume that $\{1,\ldots,K\}$ is the union of  disjoint sets $\mathcal{D}_1$, $\mathcal{D}_2$, and $\mathcal{D}_3$ such that
		\begin{itemize}
			\item[(i)] $k\in\mathcal{D}_1$ if and only if $1/R_k\in\mathcal{M}(\Phi_{\rho_k})$ for some $\rho_k\in(0,1)$.
			\item[(ii)] $k\in\mathcal{D}_2$ if and only if there exists an $\epsilon_k>0$ such that $\E \{1/R_k^{1+\epsilon_k}\}<\infty$.
			\item[(iii)] $k\in\mathcal{D}_3$ if and only if $1/R_k\in\mathcal{M}(\Phi_1)$ and $\E(1/R_k)=\infty$.
		\end{itemize}
	\end{assumption}
We conjecture that the variables whose distortions are in $\mathcal{D}_3$ have the same asymptotic behavior as those whose distortions are in $\mathcal{D}_2$. More precisely, we surmise that the following statement holds.
\begin{conjecture}\label{conj:MaxDofA}
	Let $C_{\mathcal{G},\bm{\psi},\bm{\ell},Q}$ be a clustered Archimax copula such that Assumptions~\ref{assumption:MDAExt} and \ref{assumption:stdf} hold. For $k\in\mathcal{D}_1$, let $b_k= \E(1/Z_k^{\rho_k})$, $Z_k\sim \mathrm{Beta}(1,d_k-1)$. Then $1/\bm{X} \in \mathcal{M}(H)$ with $1/X_{ki} \in \mathcal{M}(H_{ki})$, where $H_{ki}=\Phi_{\rho_k}$ for $k\in\mathcal{D}_1$ and $i \in \{1,\ldots, d_k\}$ and $H_{ki}=\Phi_{1}$ for $k\in\mathcal{D}_2 \cup \mathcal{D}_3$ and $i\in\{1,\ldots,d_k\}$. The stdf of $H$ is given for all $\bm{x}\in\mathbb{R}_+^d$ by
	\begin{align}\label{eq:stdfconj}
	\ell_{\mathcal{G},\bm{\psi},\bm{\ell},Q} (\bm{x}) =
	\E\left(\max_{\substack{ k\in\mathcal{D}_1 \\ 1 \leq i \leq d_k}}\frac{x_{ki} W_k}{b_kS_{ki}^{\rho_k}}\right)
	+ \sum_{k\in\mathcal{D}_2 \cup \mathcal{D}_3} \ell_k(\bm{x}_k)\;.
	\end{align}
\end{conjecture}

One part of Conjecture~\ref{conj:MaxDofA} is clear, namely that $H_{ki}=\Phi_{1}$ for $k \in \mathcal{D}_3$. Indeed, for any such $k$, the Corollary to Theorem~3 in \cite{Embrechts/Goldie:1980} implies that $1/(R_kS_{ki})\in\mathcal{M}(\Phi_1)$.
So one can again find  a sequence $\{a_{nk}\}$ of positive constants ensuring that for all $x\in\mathbb{R}_+$, $n\Pr\{1/(R_kS_{ki}) > a_{nk}x\} \rightarrow 1/x$ as $n\to\infty$.
The main difficulty in establishing the validity of Conjecture~\ref{conj:MaxDofA} is the fact that, for $k\in\mathcal{D}_3$ and $i\in\{1,\ldots,d_k \}$, the relation between the above normalizing sequence $\{a_{nk}\}$ and the normalizing sequences for $1/R_k$ and $1/S_{ki}$ is unclear.
In order to prove the conjectured result, it suffices to prove three sister lemmas which are reported below.
The first two, Lemmas~\ref{lemma:D1D3} and \ref{lemma:D2D3}, are analogous to Lemma~\ref{lemma:D1D2} and are proved therein.
The third, Conjecture~\ref{conjecture:D3D3}, that states asymptotic independence between different clusters in $\mathcal{D}_3$, is the missing result that if established would prove Conjecture~\ref{conj:MaxDofA}.
\begin{lemma}\label{lemma:D1D3}
	Under the hypothesis of Conjecture~\ref{conj:MaxDofA}, suppose that $k\in\mathcal{D}_1$, $l\in\mathcal{D}_3$, $i\in \{1,\ldots,d_k\}$ and $j\in\{1,\ldots,d_l \}$. Let  $\{a_{nk}\}$ be a sequence of positive constants such that for all $x>0$, $n\Pr(1/R_k>a_{nk}x)\rightarrow x^{-\rho_k}$ as $n\to\infty$ and $n\Pr\{1/(R_kS_{ki}) > a_{nk}b_k^{1/\rho_k}x\} \rightarrow x^{-\rho_k}$ as $n\to\infty$. Furthermore, let  $\{a_{nl}\}$ be  a sequence of positive constants so that for all $x > 0$, $n\Pr(1/(R_l S_{lj}) > a_{nl}x)\rightarrow 1/x$ as $n\to\infty$.
	Then for all $x,y \in\mathbb{R}_+$,
	$$
	\lim_{n\to\infty} n\Pr\{1/(R_kS_{ki}) > a_{nk}b_k^{1/\rho_k}x, 1/(R_lS_{lj})>a_{nl}y\}=0\;.
	$$
\end{lemma}

\begin{proof}[Proof of Lemma~\ref{lemma:D1D3}]
	The proof is quite similar to the one of Lemma~\ref{lemma:D1D2}. {Observe first that the assumed sequences $\{a_{nk}\}$ and $\{a_{nl}\}$ indeed exist, by Lemma \ref{lemma:marginal} and the discussion in the paragraph following Conjecture~\ref{conj:MaxDofA}. Fix some arbitrary $x,y > 0$ and} recall that $\rho_k\in(0,1)$. The probability of interest can be written as follows
	\begin{multline*}
	n\Pr\{1/(R_kS_{ki})>a_{nk}b_k^{1/\rho_k}x,\ 1/(R_lS_{lj}) > a_{nl}b_ly\}\\
	=\int_{(0,1)^2} n\Pr\{1/R_k > a_{nk} b_k^{1/\rho_k}x s_{ki},\ 1/R_l > a_{nl}y s_{lj} \} dF_{S_{ki},S_{lj}}(s_{ki},s_{lj})\;.
	\end{multline*}
	Consider the integrand as {a function $f_n$ defined on $(0,1)^2$ and note that for all $n \in \mathbb{N}$, $0 \le f_n \le g_n$, where $g_n$ is given, for all $(s_{ki},s_{lj}) \in (0,1)^2$ by}
	$$
	g_n(s_{ki},s_{lj}) = g_n(s_{ki}) = n\Pr( 1/R_k>a_{nk}b_k^{1/\rho_k}xs_{ki})\;.
	$$
	{As in the proof of Lemma~\ref{lemma:D1D2},} for all $(s_{ki},s_{lj})\in(0,1)^2$,
	$$
	\lim_{n\to\infty} g_n(s_{ki},s_{lj})= g(s_{ki},s_{lj}) = 1/\{b_k(xs_{ki})^{\rho_k}\}\;.
	$$ 
	Moreover, 
	$$
	\int_{(0,1)^2}g(s_{ki},s_{lj})dF_{S_{ki},S_{lj}}(s_{ki},s_{lj}) = \frac{1}{x^{\rho_k}}
	$$
	and
	\begin{equation*}
	\int_{(0,1)^2}g_n(s_{ki},s_{lj}) dF_{S_{ki},S_{lj}}(s_{ki},s_{lj})=n\Pr\{ 1/(R_kS_{ki})>a_{nk}b_k^{1/\rho_k}x\}\to \frac{1}{x^{\rho_k}}
	\end{equation*}
	as $n\to\infty$.
	We therefore have a sequence of functions $\{g_n\}$ bounding $\{f_n\}$ from above such that 
	\begin{align*}
		\lim_{n\to\infty}\int_{(0,1)^2} & g_n(s_{ki},s_{lj}) dF_{S_{ki},S_{lj}}(s_{ki},s_{lj})\\
		&= \int_{(0,1)^2}\lim_{n\to\infty} g_n(s_{ki},s_{lj})dF_{S_{ki},S_{lj}}(s_{ki},s_{lj})\;.
	\end{align*}
	
	Finally, note that
	\begin{equation*}
	f_n(s_{ki},s_{lj}) = n\Pr\{1/R_k > a_{nk}b_k^{1/\rho_k}x s_{ki}, 1/R_l > a_{nl}b_l^{1/\rho_l}ys_{lj}\}\to 0
	\end{equation*}
	as $n\to\infty$ since $n \Pr\{1/R_k>a_{nk}b_k^{1/\rho_k}xs_{ki}\}\to \{b_k^{1/\rho_k}xs_{ki}\}^{-\rho_k}$ and $\Pr\{1/R_l > a_{nl}b_l^{1/\rho_l}ys_{lj}\}\to 0$ as $n \to \infty$. 
	The desired result then follows from the generalized Lebesgue dominated convergence theorem.
\end{proof}

\begin{lemma}\label{lemma:D2D3}
	Under the hypothesis of Conjecture~\ref{conj:MaxDofA}, suppose that $k\in\mathcal{D}_2$, $l\in\mathcal{D}_3$, $i\in \{1,\ldots,d_k\}$ and $j\in\{1,\ldots,d_l \}$. Let  $\{a_{nk}\}$ such that for all $x>0$, $n\Pr(1/S_{ki} > a_{nk}x)\rightarrow x^{-1}$ as $n\to\infty$ and $n\Pr\{1/(R_kS_{ki})>a_{nk}{b_k}x\}\rightarrow x^{-1}$ as $n\to\infty$. Furthermore, let  $\{a_{nl}\}$ be a sequence of positive constants so that for all $x > 0$, $n\Pr\{1/(R_lS_{lj})>a_{nl}x\} \rightarrow x^{-1}$ as $n\to\infty$.
	Then for all $x,y\in\mathbb{R}_+$,
	$$
	\lim_{n\to\infty} n\Pr\{ 1/(R_kS_{ki}) > a_{nk}b_k^{1/\rho_k}x,\ 1/(R_lS_{lj}) > a_{nl}y\}=0\;.
	$$
\end{lemma}

\begin{proof}[Proof of Lemma~\ref{lemma:D2D3}]
	This proof is almost exactly the same as the proof of Lemma~\ref{lemma:D1D2}.
	Again, the existence of the norming constants  $\{a_{nk}\}$ and $\{a_{nl}\}$ follows from Lemma \ref{lemma:marginal} and the discussion in the paragraph following Conjecture~\ref{conj:MaxDofA}.
	Fix some arbitrary $x,y >0$.
	We are interested in the limit as $n\to\infty$ of 
	\begin{align*}
	& n\Pr\{1/(R_kS_{ki})>a_{nk}b_kx,\ 1/(R_lS_{lj}) > a_{nl}y\}\\
	& \qquad = \int_{\mathbb{R}_+^2} n\Pr\{1/S_{ki} > a_{nk}b_kx r_k\}\Pr\{ 1/S_{lj} > a_{nl}yr_l\} dF_{R_k,R_l}(r_k,r_l)\;.
	\end{align*}
	Consider the integrand as a {function $f_n$ defined on $\mathbb{R}_+^2$. Observe that for each $n \in \mathbb{N}$, $0 \le f_n \le g_n$ where for all $(r_k, r_l) \in \mathbb{R}_+^2$},
	$$
	g_n(r_k,r_l)=g_n(r_k)=n\Pr\{ 1/S_{ki}>a_{nk}b_kxr_k\}\;.
	$$
	From the choice of $\{a_{nk}\}$, for all $(r_k,r_l)\in\mathbb{R}_+^2$,
	$$
	\lim\limits_{n\to\infty}g_n(r_k,r_l)= g(r_k,r_l)=1/(b_kxr_k)\;.
	$$
	Moreover, {since $b_k = \E(1/R_k)$},
	$$
	\int_{\mathbb{R}_+^2}g(r_k,r_l)dF_{R_k,R_l}(r_k,r_l)=\int_{\mathbb{R}_+^2}\frac{1}{b_kxr_k}dF_{R_k,R_l}(r_k,r_l)=\frac{1}{x}\;.
	$$
	and
	$$\int_{\mathbb{R}_+^2}g_n(r_k,r_l)dF_{R_k,R_l}(r_k,r_l)=n\Pr\{ 1/(R_kS_{ki})>a_{nk}b_kx\}\to \frac{1}{x}
	$$
	as $n\to\infty$. We therefore have a sequence of functions $\{g_n\}$ bounding $\{f_n\}$ from above such that 
	$$
	\lim\limits_{n\to\infty}\int_{\mathbb{R}_+^2}g_n(r_k,r_l)dF_{R_k,R_l}(r_k,r_l)=\int_{\mathbb{R}_+^2}\lim\limits_{n\to\infty}g_n(r_k,r_l)dF_{R_k,R_l}(r_k,r_l)\;.
	$$ 
	Finally, note that
	\begin{equation*}
	f_n(r_k,r_l)=n\Pr\{1/S_{ki}>a_{nk}b_kxr_k\}\Pr\{ 1/S_{lj}>a_{nl}yr_l\}\to 0
	\end{equation*}
	as $n\to\infty$ since 
	\begin{align*}
	n\Pr\{1/S_{ki}>a_{nk}b_kxr_k\}\to 1/\{b_kxr_k \} \quad\text{and} \quad
	{\Pr\{ 1/S_{lj}>a_{nl}yr_l\}\to 0}
	\end{align*}
	as $n\to\infty$.
	Using the generalized Lebesgue dominated convergence theorem concludes the proof.
\end{proof}
\begin{conjecture}\label{conjecture:D3D3} 
	Under the hypothesis of Conjecture~\ref{conj:MaxDofA}, suppose that $k,l\in\mathcal{D}_3$, $i\in \{1,\ldots,d_k\}$ and $j\in\{1,\ldots,d_l \}$. Let  $\{a_{nk}\}$ and $\{a_{nl}\}$ be sequences of positive constants such that for all $x>0$, $n\Pr(1/(R_kS_{ki})>a_{nk}x)\rightarrow x^{-1}$ and $n\Pr(1/(R_lS_{lj})>a_{nl}x)\rightarrow x^{-1}$  as $n\to\infty$.
	Then for all $x,y\in\mathbb{R}_+$,
	$$
	\lim_{n\to\infty} n\Pr\{1/(R_kS_{ki}) > a_{nk}b_k^{1/\rho_k}x,\ 1/(R_lS_{lj})>a_{nl}y\}=0\;.
	$$
\end{conjecture}

\section{Additional figures for Section~\ref{sec:simulations}}\label{sec:additionalfigures}
\begin{figure}[h!]
	\includegraphics[scale=0.4]{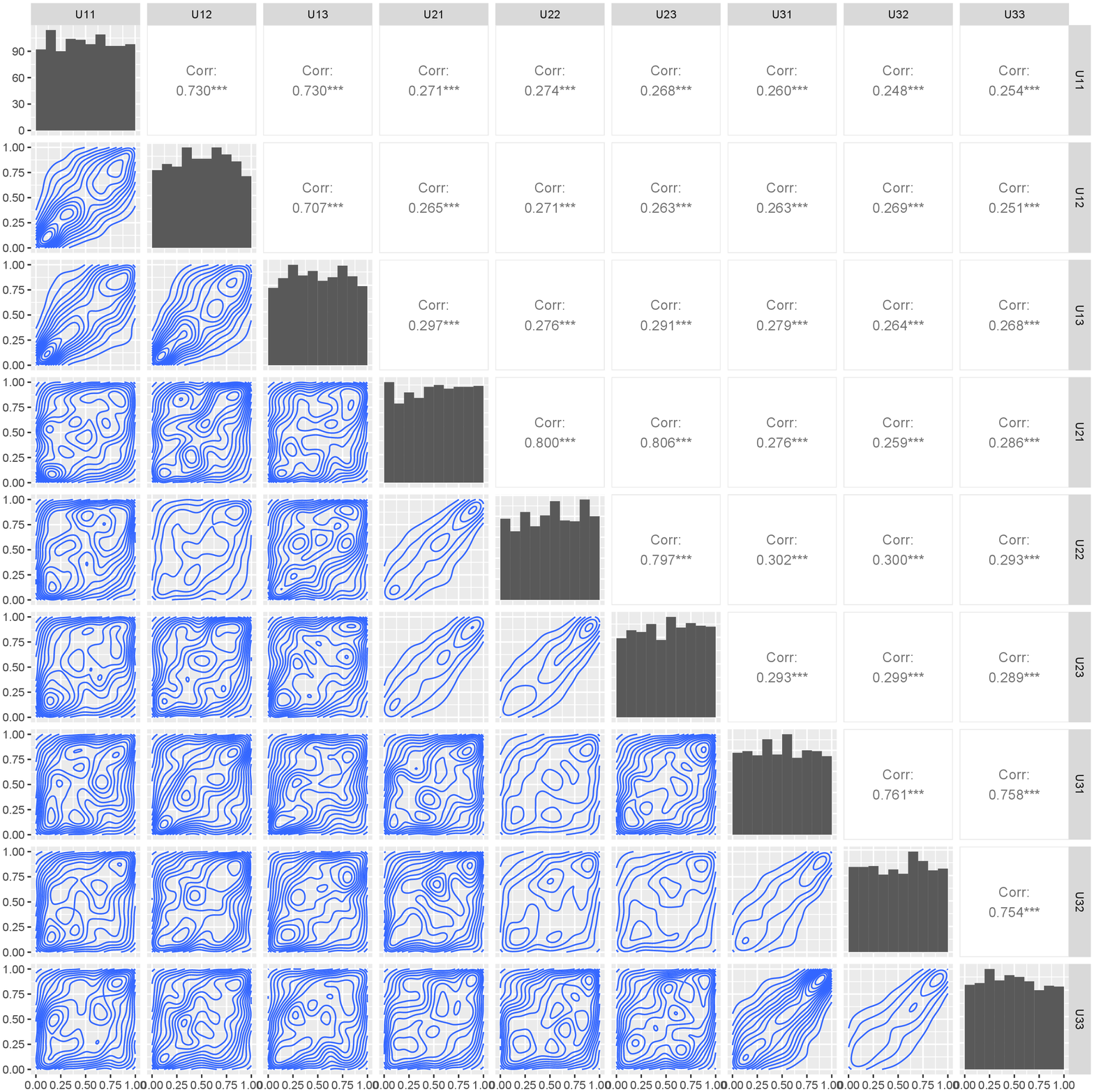}
	\caption{(Model A) Pairwise plots of a sample of size $n=1000$ from the copula $C_{\mathcal{G},\bm{\psi},\bm{\ell},Q}$ where  $\mathcal{G}=\{\mathcal{G}_1,\mathcal{G}_2,\mathcal{G}_3 \}=\{\{1,2,3\},\{4,5,6\},\{7,8,9\} \}$. $\bm{\psi}=\{\psi_{\theta_1},\psi_{\theta_2},\psi_{\theta_3} \}$ with $(\theta_1,\theta_2,\theta_3)=(1.5,1.5,2)$ where $\psi_{\theta_1}$ is Clayton while $\psi_{\theta_2},\psi_{\theta_3}$ are Joe. $\bm{\ell}=\{\ell_{\vartheta_1},\ell_{\vartheta_2},\ell_{\vartheta_3} \}$ with $(\vartheta_1,\vartheta_2,\vartheta_3)=(1.25,2,1.5)$ where $\ell_{\vartheta_1},\ell_{\vartheta_2},\ell_{\vartheta_3}$ are Gumbel-Hougaard. The radial survival copula $\bar{D}$ is trivariate Gaussian with correlations all equal to $0.5$. Upper: linear correlation, Lower: contour density, Diagonal: univariate histogram.}
	\label{fig:pairplotN}
\end{figure}

\begin{figure}[h!]
	\includegraphics[scale=0.4]{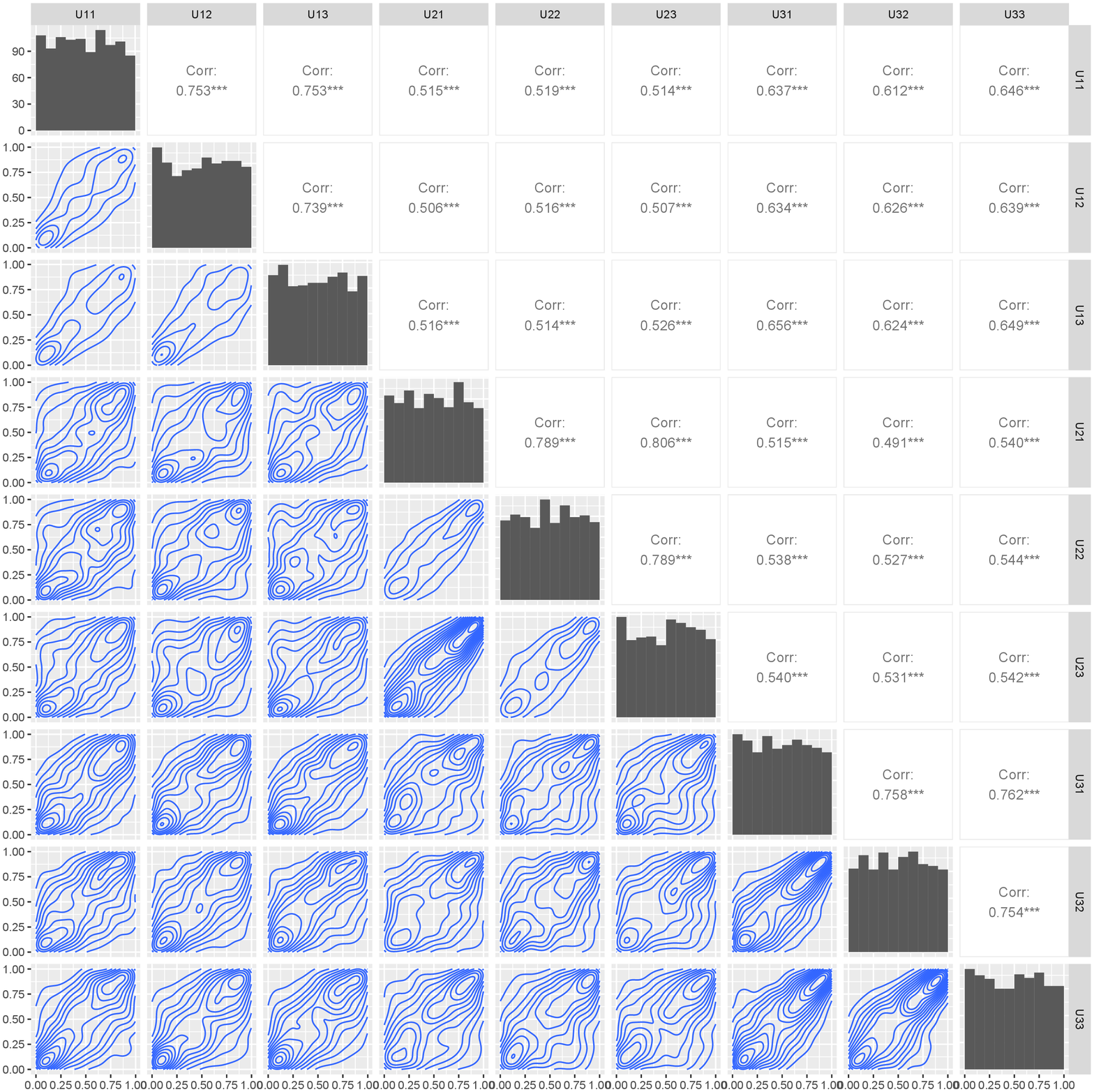}
	\caption{(Model B) Pair plots of a sample of size $n=1000$ from the copula $C_{\mathcal{G},\bm{\psi},\bm{\ell},Q}$ where  $\mathcal{G}=\{\mathcal{G}_1,\mathcal{G}_2,\mathcal{G}_3 \}=\{\{1,2,3\},\{4,5,6\},\{7,8,9\} \}$. $\bm{\psi}=\{\psi_{\theta_1},\psi_{\theta_2},\psi_{\theta_3} \}$ with $(\theta_1,\theta_2,\theta_3)=(1.5,1.5,2)$ where $\psi_{\theta_1}$ is Clayton while $\psi_{\theta_2},\psi_{\theta_3}$ are Joe. $\bm{\ell}=\{\ell_{\vartheta_1},\ell_{\vartheta_2},\ell_{\vartheta_3} \}$ with $(\vartheta_1,\vartheta_2,\vartheta_3)=(1.25,2,1.5)$ where $\ell_{\vartheta_1},\ell_{\vartheta_2},\ell_{\vartheta_3}$ are Gumbel-Hougaard. The radial survival copula $\bar{D}$ is trivariate Gumbel-Hougaard with parameter $\vartheta_R=4$. Upper: linear correlation, Lower: contour density, Diagonal: univariate histogram.}
	\label{fig:pairplotG}
\end{figure}
\begin{figure}
	\begin{subfigure}[b]{0.45\linewidth}
		\includegraphics[scale=0.5]{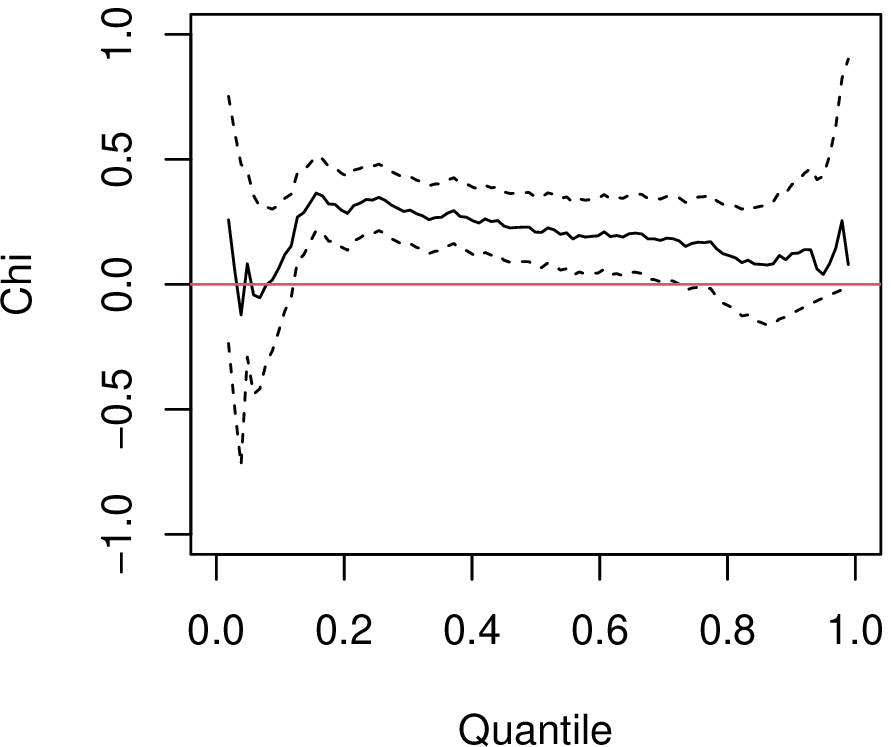}
		\caption{Model A: Pair $\{1,4\}$}
	\end{subfigure}
	\begin{subfigure}[b]{0.45\linewidth}
		\includegraphics[scale=0.5]{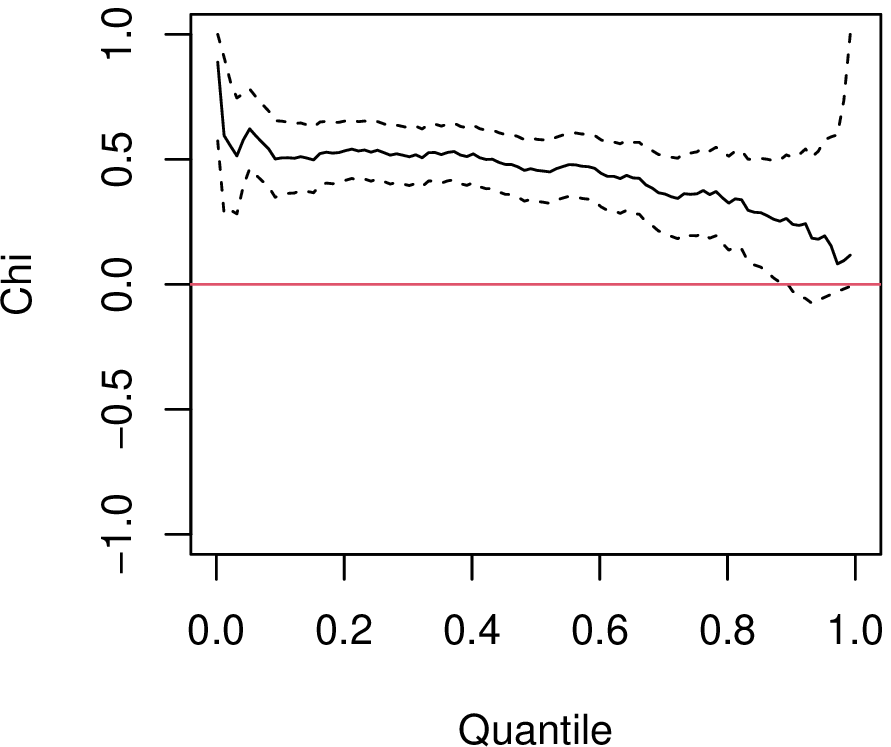}
		\caption{Model B: Pair $\{1,4\}$}
	\end{subfigure}
	\begin{subfigure}[b]{0.45\linewidth}
		\includegraphics[scale=0.5]{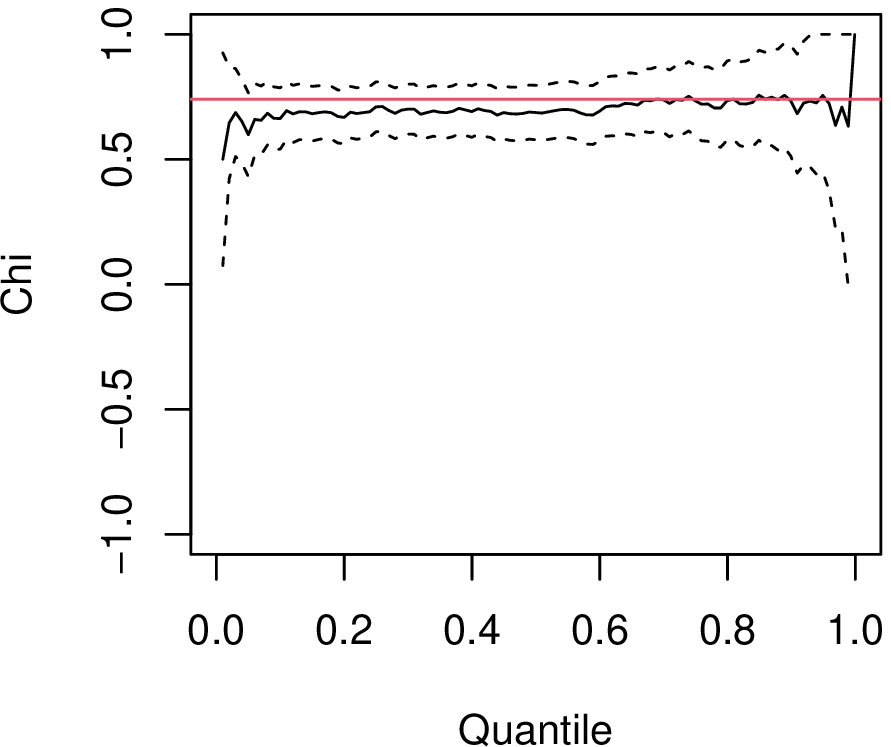}
		\caption{Model A: Pair $\{4,5\}$}
	\end{subfigure}
	\begin{subfigure}[b]{0.45\linewidth}
		\includegraphics[scale=0.5]{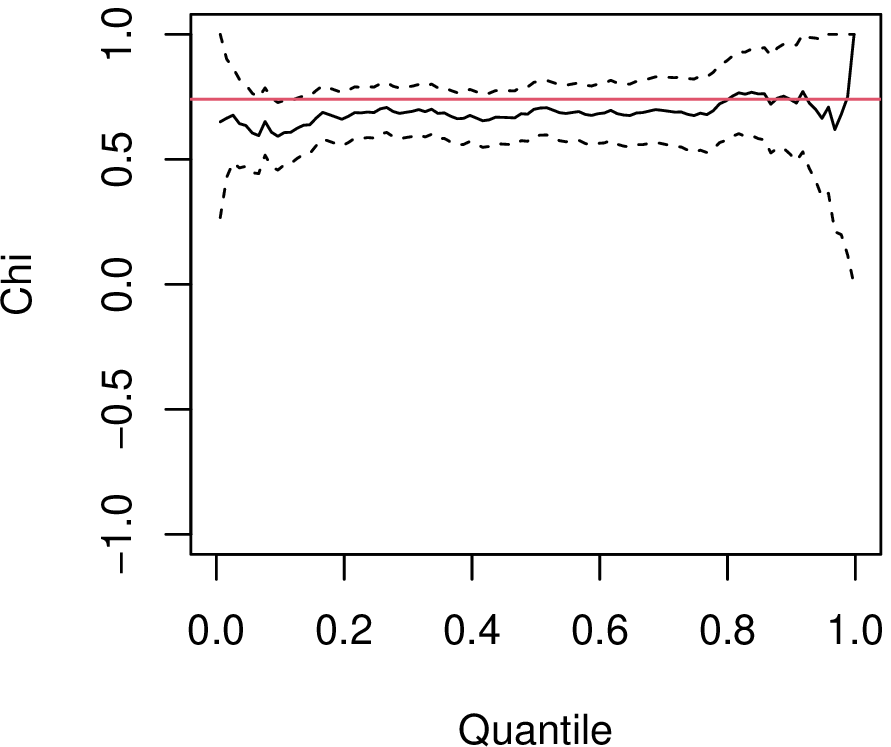}
		\caption{Model B: Pair $\{4,5\}$}
	\end{subfigure}
	\begin{subfigure}[b]{0.45\linewidth}
		\includegraphics[scale=0.5]{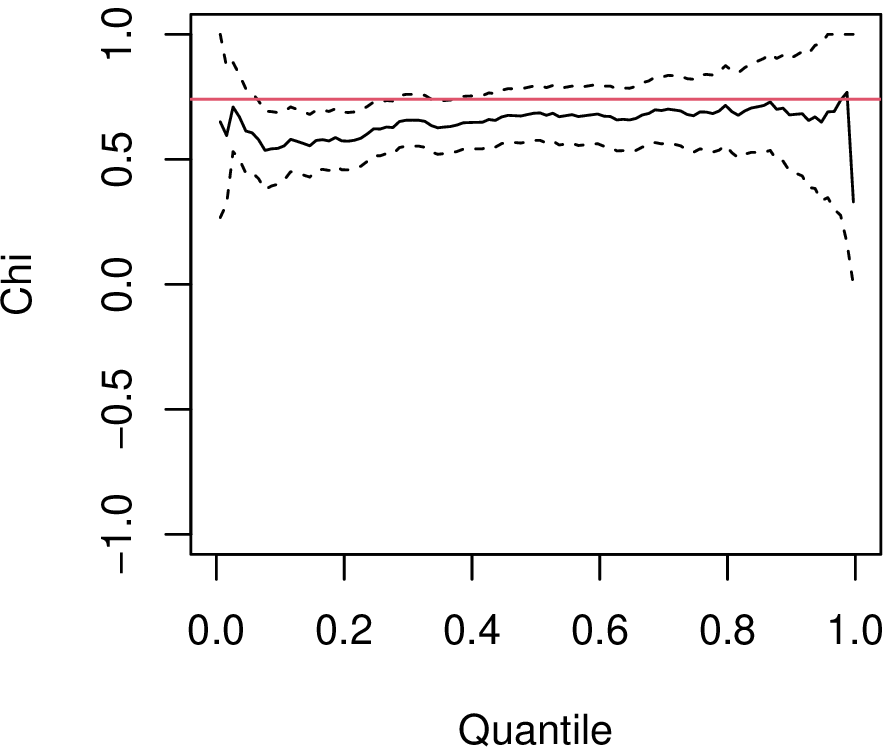}
		\caption{Model A: Pair $\{7,8\}$}
	\end{subfigure}
	\begin{subfigure}[b]{0.45\linewidth}
		\includegraphics[scale=0.5]{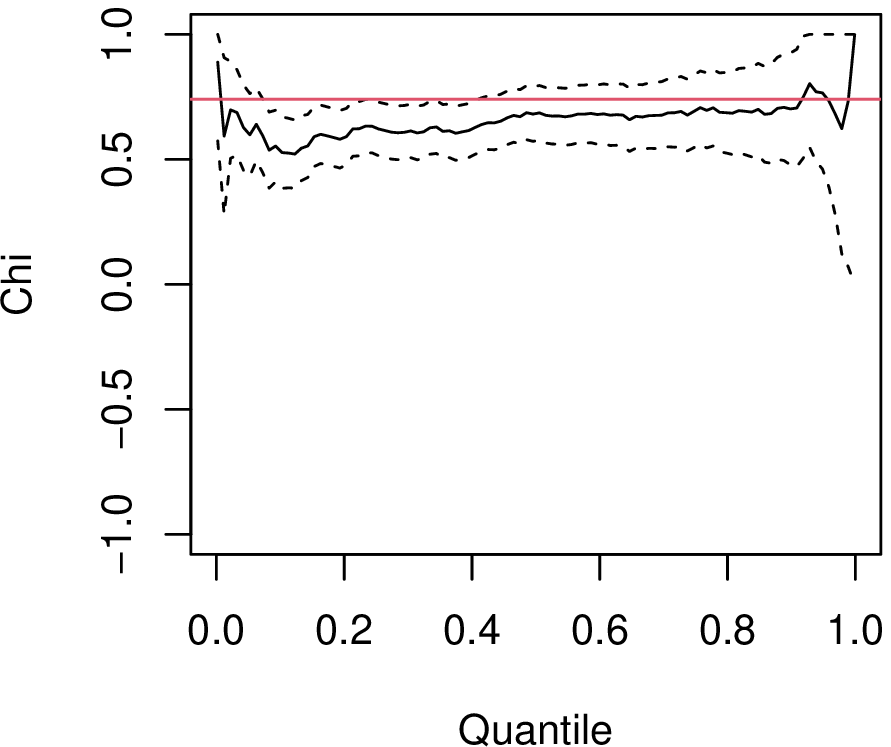}
		\caption{Model B: Pair $\{7,8\}$}
	\end{subfigure}
	\caption{ Pair chi plots for the sample from Model A (left) and Model B (right). The full black line is the empirical estimate of \eqref{eq:chi}, the dotted lines are 95\% confidence intervals and the red lines represent the true values of $\lim_{q\to 1}\chi_{ij}(q)$ for each pair ${(i,j)}$. The samples used for the empirical estimates are those of Figures~\ref{fig:pairplotN} and \ref{fig:pairplotG}}
	\label{fig:chiplots2}
\end{figure}
\clearpage
\section{Additional figures for the data illustration of Section~\ref{sec:data-application}}

\bigskip
\bigskip

\begin{figure}[h]
	\includegraphics[width=1\linewidth]{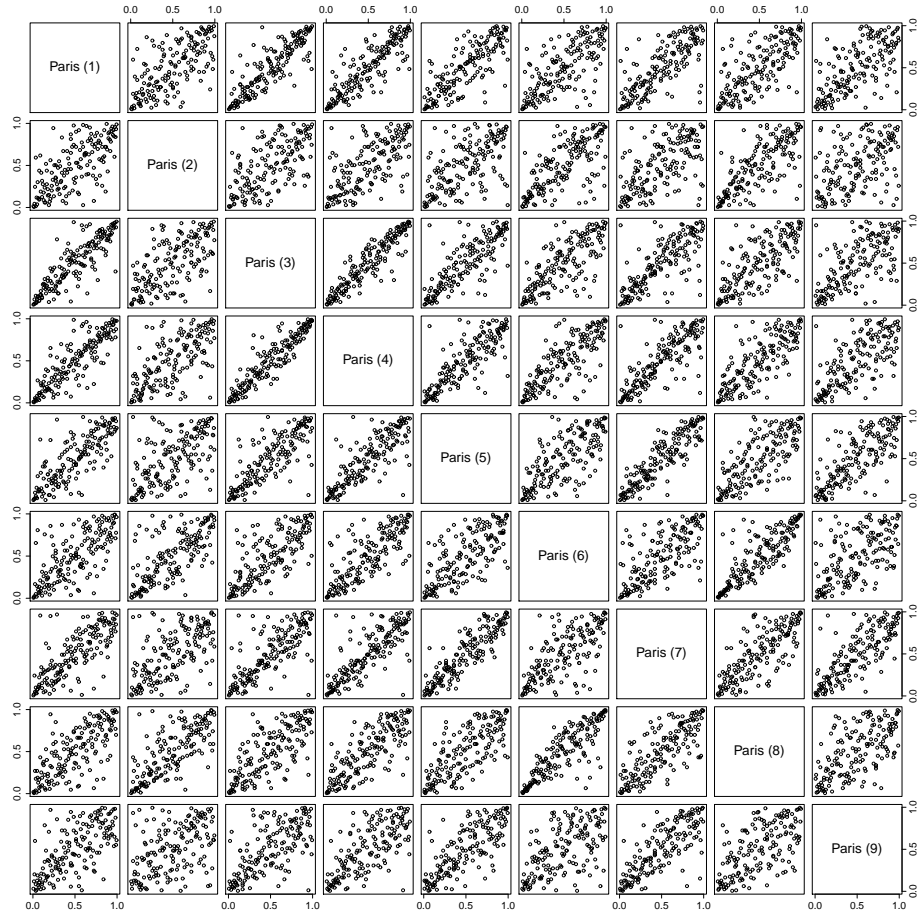}
	\caption{Pair plots of the scaled componentwise ranks of monthly maxima of precipitations for the nine stations of the Paris cluster considered in the application of Section~\ref{sec:data-application}.} \label{fig:paris}
\end{figure}
\begin{figure}[h]
	\includegraphics[width=1\linewidth]{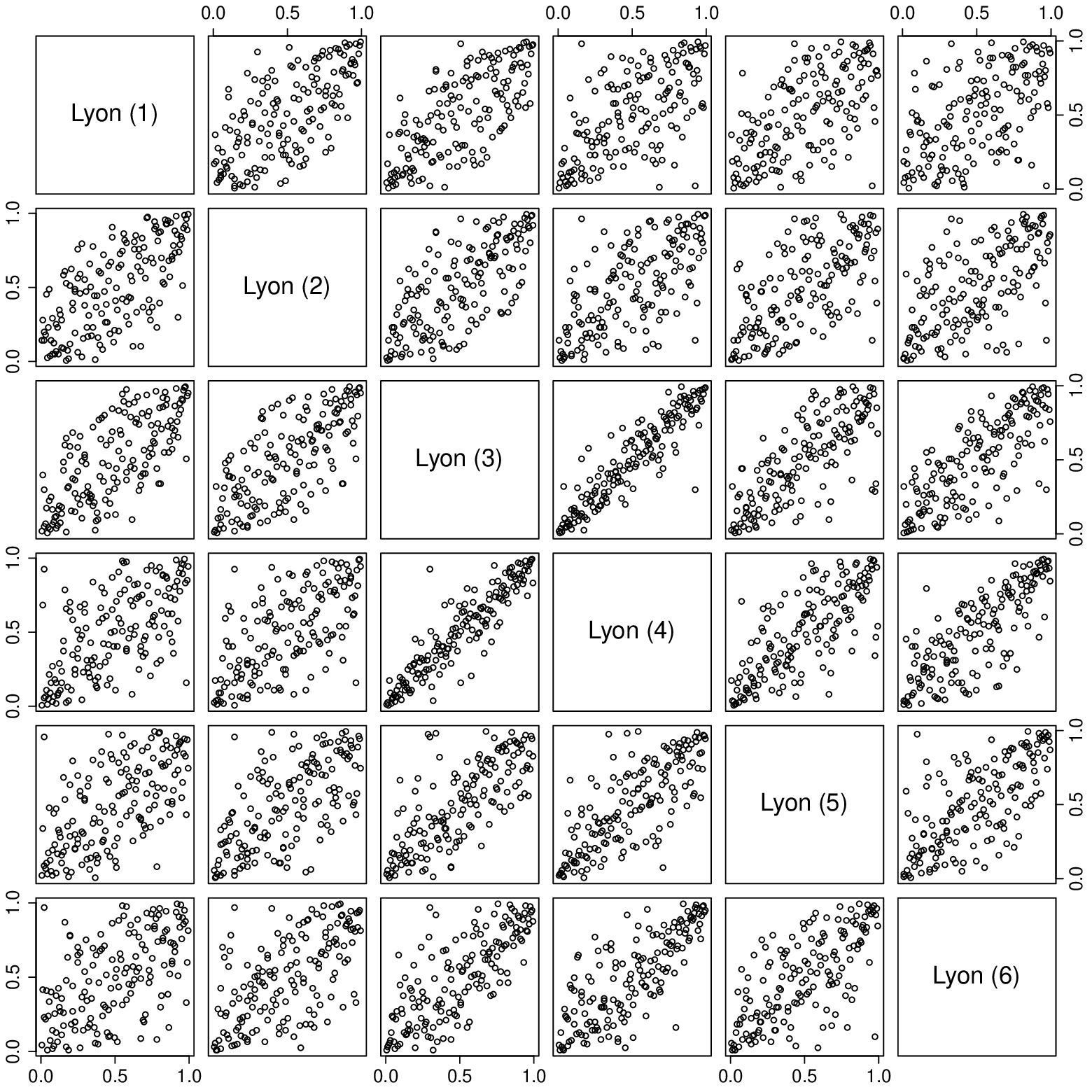}
	\caption{Pair plots of the scaled componentwise ranks of monthly maxima of precipitations for the six stations of the Paris cluster considered in the application of Section~\ref{sec:data-application}.} \label{fig:lyon}
\end{figure}
\begin{figure}[h]
	\includegraphics[width=1\linewidth]{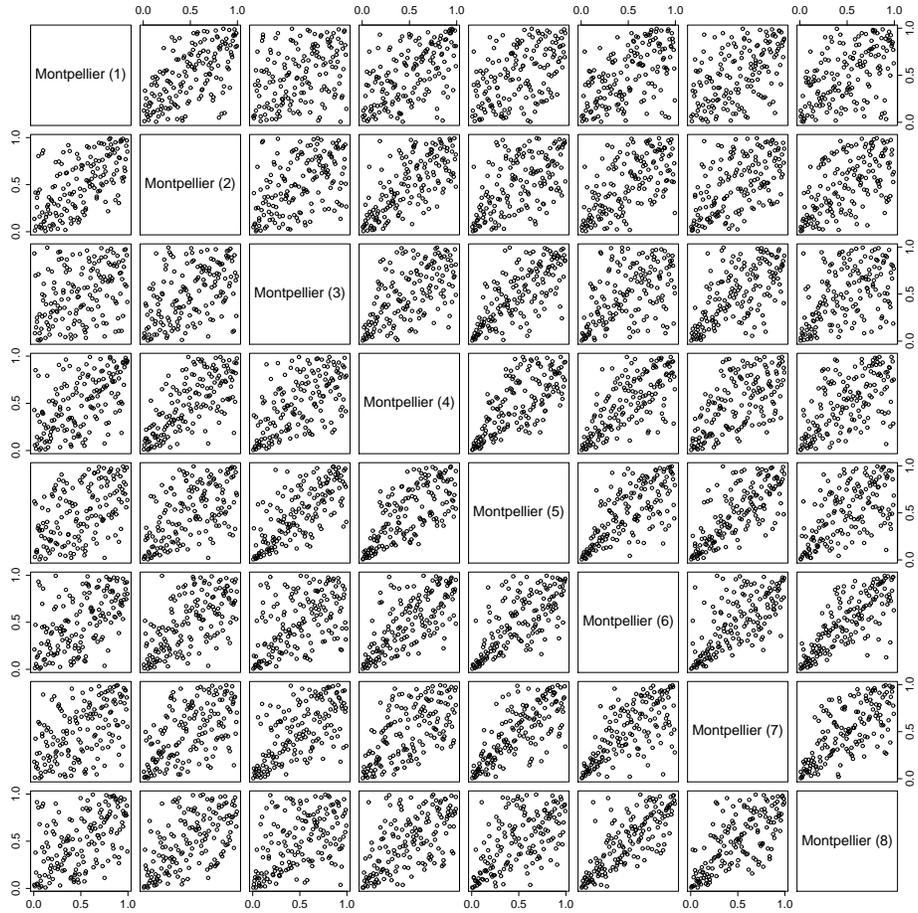}
	\caption{Pair plots of the scaled componentwise ranks of monthly maxima of precipitations for the eight original stations of the Montpellier cluster considered in the application of Section~\ref{sec:data-application}.} \label{fig:montpellier}
\end{figure}

\begin{figure}[t]
\centering
\includegraphics[width=.48\linewidth]{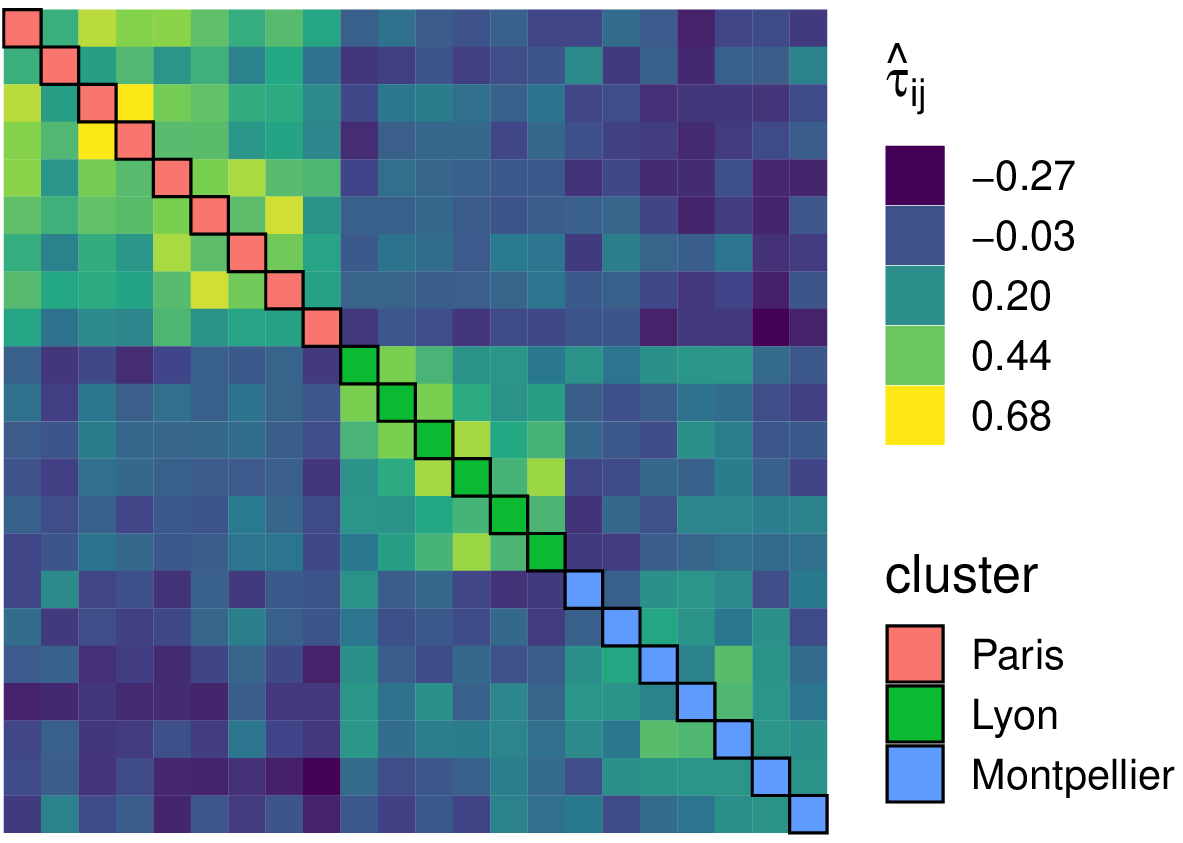}
\quad
\includegraphics[width=.48\linewidth]{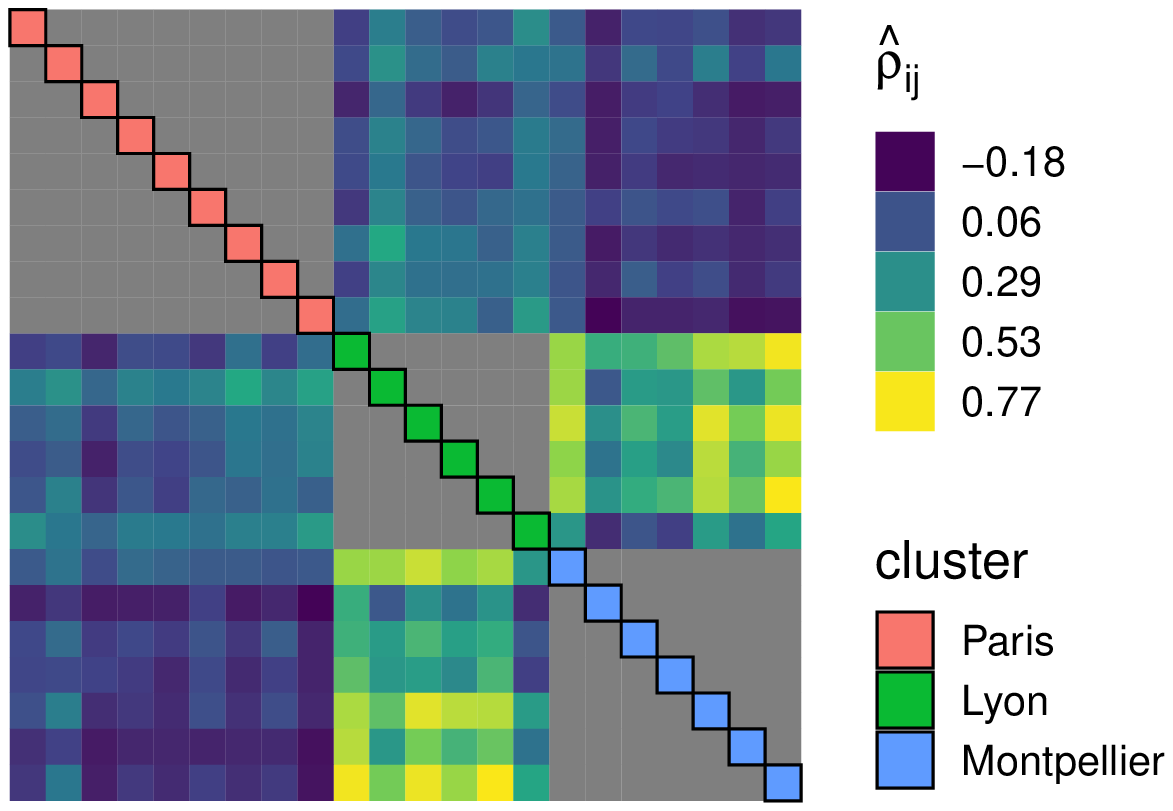}
\caption{Estimated quantities from Section~\ref{sec:app-inference-inter} ($d=22$ stations). Left: Matrix of pairwise Kendall correlations associated with the (yearly) data used to perform the test of \cite{Kojadinovic/Holmes:2009}. Right: Matrix of pairwise Pearson correlations used for estimating $Q$.} \label{fig:rho-tau-hat}
\end{figure}

\end{appendices}
\end{document}